\documentclass{article}
\usepackage{amsmath}
\usepackage{amsmath,amssymb,latexsym,color}
\usepackage[mathscr]{eucal}
\newtheorem{thm}{Theorem}[section]

\newtheorem{lemma}[thm]{Lemma}

\newtheorem{example}{Example}[section]
\newtheorem{defin}{Definition}[section]

\newtheorem{remark}{Remark}[section]

\newcommand{\proof}{{\it Proof.\quad}}
\newcommand{\qed}{\hfill\Box\medskip}

\usepackage{CJK}
\begin{document}
\begin{CJK*}{GBK}{song}

\renewcommand{\baselinestretch}{1.2}

\title{\bf Cameron-Liebler sets in bilinear forms graphs}

\author{
Jun Guo\thanks{Corresponding author. guojun$_-$lf@163.com}\\
{\footnotesize Department of Mathematics, Langfang Normal University, Langfang  065000,  China} }
 \date{ }
 \maketitle

\begin{abstract}
  Cameron-Liebler sets of subspaces in projective spaces were studied recently by
Blokhuis, De Boeck and D'haeseleer (Des. Codes Cryptogr., 2019). In this paper,
  we discuss Cameron-Liebler sets in bilinear forms graphs, obtain several equivalent
definitions and present some classification results.



\end{abstract}

\section{Introduction}
\underline{}Cameron-Liebler sets of lines were first introduced by Cameron and Liebler \cite{Cameron} in their study of collineation groups of PG$(3,q)$. There have been many results for Cameron-Liebler sets of lines in the projective space PG$(3,q)$. See \cite{Metsch2,Metsch3,Rodgers} for classification results,  and \cite{Bruen,De Beule,Feng,Gavrilyuk} for the constructions of two non-trivial examples.
Over the years, there have been many interesting extensions of this notion.
See \cite{De Boeck2} for Cameron-Liebler classes in finite sets, \cite{De Boeck} for Cameron-Liebler sets of generators in polar spaces, \cite{Blokhuis,Gavrilyuk2,Metsch,Rodgers2}  for Cameron-Liebler sets of $k$-spaces in PG$(n,q)$.

Similar problems have been investigated by various researchers under different names. See \cite{Martin} for completely regular designs of strength 0, \cite{Bamberg} for tight sets, \cite{De Bruyn} for intriguing sets, and \cite{Filmus} for boolean degree 1 functions. One of the main reasons for studying Cameron-Liebler sets is that there are many connections to other geometric and combinatorial objects, such as blocking sets, intersecting families, linear codes, and association schemes.

For nonnegative integers $n$ and $\ell$, let
$\mathbb{F}_q^{n+\ell}$ be an $(n+\ell)$-dimensional column vector space over the finite field $\mathbb{F}_q$.
Let $e_i\,(1\leq i\leq n+\ell)$ be the column vector in $\mathbb{F}_q^{n+\ell}$ whose $i$th coordinate is 1 and all other coordinates are 0. Denote by $E=\langle e_{n+1},\ldots,e_{n+\ell}\rangle$ the $\ell$-dimensional subspace spanned by $e_{n+1},\ldots,e_{n+\ell}$. An $m$-dimensional subspace $P$ of $\mathbb{F}_q^{n+\ell}$ is called a subspace of type $(m,k)$ if $\dim(P\cap E)=k$. Denote the set of all subspaces of type $(m,k)$ in $\mathbb{F}_q^{n+\ell}$ by ${\cal M}(m,k;n+\ell,E)$. The {\it attenuated space} $A_q(n+\ell,E)$ is the collection of all subspaces of $\mathbb{F}_q^{n+\ell}$ intersecting trivially with $E$.

For convenience, we write ${\cal M}(m,0;n+\ell,E)$ as ${\cal M}_m$ for $m=1,2,\ldots,n$.
Then $A_q(n+\ell,E)=\bigcup_{i=0}^n{\cal M}_i$.
Let $1\leq i\leq j\leq n$. For a fixed subspace $\tau\in{\cal M}_i$, let ${\cal M}_j'(\tau)$ be the set
of all subspaces in  ${\cal M}_j$ containing $\tau$.
For a fixed subspace $\pi\in{\cal M}_n$, let
$\overline{{\cal M}}_i(\pi)$ be the set of all subspaces in ${\cal M}_i$ disjoint to $\pi$.
For a fixed subspace $\Sigma$ of $\mathbb{F}_q^{n+\ell}$, let ${\cal M}_i(\Sigma)$ be the set of all subspaces in  ${\cal M}_i$ contained in $\Sigma$.

In 2010, Wang, Guo and Li \cite{Wang} proved that the set ${\cal M}_m$ has a structure of a symmetric association scheme $\mathfrak{X}_{m,n,\ell}$
and computed intersection numbers of $\mathfrak{X}_{m,n,\ell}$. The scheme $\mathfrak{X}_{m,n,\ell}$ is a common generalization of
the Grassmann scheme $\mathfrak{X}_{m,n,0}$ and the bilinear forms scheme $\mathfrak{X}_{n,n,\ell}$. In 2019, Blokhuis, De Boeck and D'haeseleer
\cite{Blokhuis} studied Cameron-Liebler sets in the Grassmann scheme $\mathfrak{X}_{m,n,0}$. Their  research stimulates us to consider
 Cameron-Liebler sets in the bilinear forms scheme $\mathfrak{X}_{n,n,\ell}$.

The {\it bilinear forms graph} ${\rm Bil}_q(n,\ell)$ has the vertex set
${\cal M}_n$, and two vertices are adjacent if their intersection is of dimension $n-1$.
The bilinear forms graph ${\rm Bil}_q(n,\ell)$ is
a distance-regular graph with $q^{n\ell}$ vertices and diameter $\min\{n,\ell\}$.
Note that the distance between two vertices $A$ and $B$ in ${\cal M}_n$ is $n-\dim(A\cap B)$.
Since ${\rm Bil}_q(n,\ell)$ is isomorphic to ${\rm Bil}_q(\ell,n)$, we always assume that $n\leq \ell$ in the rest of this paper.

We always assume that vectors are regarded as column vectors. For any vector $\alpha$ whose positions correspond to elements in a set, we denote its value on the position corresponding to an element $a$ by $(\alpha)_a$.
The {\it characteristic vector} $\chi_{{\cal S}}$ of a subset ${\cal S}$ of ${\cal M}_n$ is the vector whose positions correspond to the elements of ${\cal M}_n$, such that $(\chi_{{\cal S}})_{\pi}=1$ if $\pi\in{\cal S}$ and 0 otherwise. The all-one vector will be denoted by ${\rm j}$.
Let $M$ be the incidence matrix of one-dimensional subspaces versus $n$-dimensional subspaces of $A_q(n+\ell,E)$.

 A subset ${\cal L}$ of ${\cal M}_n$ is called a {\it Cameron-Liebler set} in the bilinear forms graph ${\rm Bil}_q(n,\ell)$ with parameter $x=q^{-(n-1)\ell}|{\cal L}|$ if $\chi_{{\cal L}}\in {\rm Im}(M^t)$,
 where $M^t$ is the transpose matrix of $M$. A family ${\cal F}\subseteq{\cal M}_n$ is called {\it intersecting} if $\dim(\pi\cap\pi')\geq 1$ for all $\pi,\pi'\in{\cal F}$. A Cameron-Liebler set ${\cal L}$ in ${\rm Bil}_q(n,\ell)$ with parameter $x$ is called {\it trivial} if ${\cal L}$ is a union of $x$ intersecting families.

 Cameron-Liebler sets in bilinear forms graphs are studied recently.
 Bailey, Cameron, Gavrilyuk and Goryainov \cite{Bailey} studied Cameron-Liebler sets  in ${\rm Bil}_2(2,\ell)$. Filmus and Ihringer \cite{Filmus} made a  guess on how all Cameron-Liebler sets in ${\rm Bil}_q(n,\ell)$ should look like and obtained a non-trivial example for $n=\ell=2$.
Ihringer \cite{Ihringer} obtained some classification results for Cameron-Liebler sets in ${\rm Bil}_q(n,\ell)$.

In this paper, we consider Cameron-Liebler sets in the bilinear forms graph ${\rm Bil}_q(n,\ell)$.
The rest of this paper is structured as follows.
In Section~2, we give several equivalent
definitions for these Cameron-Liebler sets in ${\rm Bil}_q(n,\ell)$.
In Section~3, we obtain several examples and properties of these Cameron-Liebler sets in ${\rm Bil}_q(n,\ell)$.
By using the properties from Section~3, we give the following stronger classification result in Section~4:
 If $\ell\geq2n\geq4$ and $2\leq x\leq (q-1)^{\frac{n}{2}}q^{\frac{\ell-2n+1}{2}}$, then there are no non-trivial Cameron-Liebler sets in ${\rm Bil}_q(n,\ell)$ with parameter $x$.

\section{Several equivalent definitions}
Let $G$ be the graph  with the vertex set ${\cal M}_1$, and two vertices $\tau$ and $\tau'$ are adjacent if $\tau+\tau'\in{\cal M}_2$. Then $G$ is a complete ${n\brack 1}_q$-partite graph.

\begin{lemma}\label{lem2.2}{\rm(See Lemma~2 in \cite{Esser}.)}
The distinct eigenvalues of $G$ are $q^{\ell+1}{n-1\brack 1}_q,0$ and $-q^{\ell}$,
and the corresponding multiplicities are $1,(q^{\ell}-1){n\brack 1}_q$ and
$q{n-1\brack 1}_q$, respectively.
\end{lemma}

\begin{lemma}\label{lem2.1}{\rm(See Theorem~3.3 in \cite{Hirschfeld}.)}
Let $1\leq i\leq j\leq n$ and $\tau\in{\cal M}_i$. Then the size of the set ${\cal M}'_j(\tau)$ is
$q^{\ell(j-i)}{n-i\brack j-i}_q.$
\end{lemma}

\begin{lemma}\label{lem2.3}
The rank of the incidence matrix $M$ is  $(q^{\ell}-1){n\brack 1}_q+1$ over the real field $\mathbb{R}$.
\end{lemma}
\proof Let $N=MM^t$. Then both $N$ and
$M$ have the same rank over $\mathbb{R}$. Note that $N$ is a $q^{\ell}{n\brack 1}_q\times q^{\ell}{n\brack 1}_q$ matrix with rows and columns indexed with the subspaces in ${\cal M}_1$.
 For any $\tau,\tau'\in{\cal M}_1$, the entry
$N_{\tau\tau'}$ is the number of subspaces of type $(n,0)$ containing both $\tau$ and $\tau'$. By Lemma~\ref{lem2.1}, we have
$$N_{\tau\tau'}=\left\{\begin{array}{ll}
q^{(n-1)\ell}  &  \hbox{if}\;\tau=\tau',\\
q^{(n-2)\ell}   & \hbox{if}\;\tau+\tau'\in{\cal M}_2,\\
0           &   \hbox{otherwise},
\end{array}\right.$$
which implies that
\begin{equation}\label{equa1}
N=q^{(n-1)\ell}I+q^{(n-2)\ell}A,
\end{equation}
where $I$ is the identity matrix of order $q^{\ell}{n\brack 1}_q$,
and $A$ is the adjacency matrix of the graph $G$. By Lemma~\ref{lem2.2} and (\ref{equa1}), we obtain that
the distinct eigenvalues of $N$ are $q^{(n-1)\ell}{n\brack 1}_q,q^{(n-1)\ell}$ and $0$, and
the corresponding multiplicities are $1,(q^{\ell}-1){n\brack 1}_q$ and
$q{n-1\brack 1}_q$, respectively.
It follows  that the rank of $N$ is $(q^{\ell}-1){n\brack 1}_q+1$. $\qed$

The {\it attenuated $q$-Kneser graph} $AK_q(n,\ell)$ has the vertex set
${\cal M}_n$, and two vertices are adjacent if they are disjoint.
Let $K$ be the adjacent matrix of  $AK_q(n,\ell)$.
The eigenvalues and the dimensions of the eigenspaces of $AK_q(n,\ell)$ are described in \cite{Delsarte},
these formulas are simplified in \cite{Lv}. So,
 we will use the following simplified formulas.

\begin{lemma}\label{lem2.4}{\rm(See Theorem~8 in \cite{Lv}.)}
The distinct eigenvalues of $AK_q(n,\ell)$ are
$\lambda_{j}=(-1)^jq^{{n\choose 2}}\prod_{s=1}^{n-j}(q^{\ell-n+s}-1),\;j=0,1,\ldots,n,$
and the eigenspace $V_j$ corresponding to $\lambda_j$ has dimension ${n\brack j}_q\prod_{s=0}^{j-1}(q^\ell-q^s)$.
\end{lemma}

\begin{lemma}\label{lem2.5.2}{\rm(See Theorem~1.10 in \cite{Wan}.)}
Let $0\leq m\leq \min\{n,\ell\}$. Then the number of $n\times \ell$ matrices of rank $m$ over $\mathbb{F}_q$ is
$q^{{m\choose 2}}{n\brack m}_q\prod_{s=1}^m(q^{\ell-m+s}-1).$
\end{lemma}

Let $P$ be an $m$-dimensional subspace of $\mathbb{F}_q^{n+\ell}$. Denote also by $P$ an $(n+\ell)\times m$ matrix of rank $m$ whose columns span the subspace $P$ and call the matrix $P$ a {\it matrix representation} of the subspace $P$.

\begin{lemma}\label{lem2.6}
Let $\tau\in{\cal M}_1$ and $\pi\in{\cal M}_n$. Then the following hold:
\begin{itemize}
\item[\rm(i)]
The size of $\overline{{\cal M}}_n(\pi)$ is $q^{{n\choose 2}}\prod_{s=1}^n(q^{\ell-n+s}-1)$.

\item[\rm(ii)]
The size of $\overline{{\cal M}}_n(\pi)\cap{\cal M}'_n(\tau)$ is $0$ if $\tau\subseteq\pi$, and $q^{{n\choose 2}}\prod_{s=1}^{n-1}(q^{\ell-n+s}-1)$ otherwise.
\end{itemize}
\end{lemma}
\proof
(i). By Theorem~9.5.2 in \cite{Brouwer},  the size of $\overline{{\cal M}}_n(\pi)$ is
$q^{{n\choose 2}}\prod_{s=1}^{n}(q^{\ell-n+s}-1)$.

(ii). If $\tau\subseteq\pi$, then the result is obvious. Suppose $\tau\not\subseteq\pi$. By the transitivity of the group $P\Gamma\!L(n+\ell,\mathbb{F}_q)_E$ on $\overline{{\cal M}}_n(\pi)$ (cf. \cite{Deng}), we may assume that $\tau=\langle e_1+e_{n+1}\rangle$ and $\pi=\langle e_1,\ldots,e_n\rangle$.
Then $\pi'\in\overline{{\cal M}}_n(\pi)\cap{\cal M}'_n(\tau)$ has a matrix representation
$$\bordermatrix{ &\hbox{\footnotesize{$1$}}&\hbox{\footnotesize{$n-1$}}\cr
&1&0\cr
&0&I\cr
&1&A_{1}\cr
&0&A_2}
\hspace{-3pt}
\begin{array}{c}
\hbox{\footnotesize{$1$}}\\
\hbox{\footnotesize{$n-1$}}\\
\hbox{\footnotesize{$1$}}\\
\hbox{\footnotesize{$\ell-1$}}
\end{array}\;\hbox{with rank}(A_2)=n-1.$$
By Lemma~\ref{lem2.5.2}, the size of $\overline{{\cal M}}_n(\pi)\cap{\cal M}'_n(\tau)$ is
$$q^{n-1}q^{{n-1\choose 2}}\prod_{s=1}^{n-1}(q^{\ell-n+s}-1)=q^{{n\choose 2}}\prod_{s=1}^{n-1}(q^{\ell-n+s}-1).$$
The proof is completed. $\qed$

\begin{lemma}\label{lem2.5}
${\rm Im}(M^t)=V_0\oplus V_1$, where $V_0=\langle {\rm j}\rangle$.
\end{lemma}
\proof By Lemma~\ref{lem2.3}, the matrix $M$ has ${n\brack 1}_qq^\ell$
rows with rank$(M)=(q^{\ell}-1){n\brack 1}_q+1$. By Lemma~\ref{lem2.4},
$\dim V_1=(q^{\ell}-1){n\brack 1}_q$, and therefore $\dim(V_0\oplus V_1)=\dim{\rm Im}(M^t)$.

 From Lemma~\ref{lem2.6}, we deduce that $K\chi_{{\cal M}'_n(\tau)}=-\lambda_1({\rm j}-\chi_{{\cal M}'_n(\tau)})$ and $K{\rm j}=-(q^\ell-1)\lambda_1{\rm j}$, which imply that
 $$K(\chi_{{\cal M}'_n(\tau)}-q^{-\ell}{\rm j})=\lambda_1(\chi_{{\cal M}'_n(\tau)}-q^{-\ell}{\rm j}).$$
It follows that $\chi_{{\cal M}'_n(\tau)}-q^{-\ell}{\rm j}\in V_1$. Therefore, we have $\chi_{{\cal M}'_n(\tau)}\in V_0\oplus V_1$.
Since $\chi_{{\cal M}'_n(\tau)}$ is the column of $M^t$ corresponding to the subspace $\tau$,
we have  ${\rm Im}(M^t)\subseteq V_0\oplus V_1$. From $\dim(V_0\oplus V_1)=\dim{\rm Im}(M^t)$, we deduce that ${\rm Im}(M^t)=V_0\oplus V_1$. $\qed$

The {\it incidence vector} $v_{{\cal S}}$ of a subset ${\cal S}$ of ${\cal M}_1$ is the vector whose positions correspond to the elements of ${\cal M}_1$, such that $(v_{{\cal S}})_{\tau}=1$ if $\tau\in{\cal S}$ and 0 otherwise.

\begin{lemma}\label{lem2.7}
Let $\pi\in{\cal M}_n$.
Then $$\chi_{\overline{{\cal M}}_n(\pi)}-(q^{-(n-1)\ell}{\rm j}-\chi_{\{\pi\}})q^{{n\choose 2}}\prod_{s=1}^{n-1}(q^{\ell-n+s}-1)\in\ker(M).$$
\end{lemma}
\proof By Lemma~\ref{lem2.6}, we have
$$M\chi_{\overline{{\cal M}}_n(\pi)}=({\rm j}-v_{{\cal M}_1(\pi)})q^{{n\choose 2}}\prod_{s=1}^{n-1}(q^{\ell-n+s}-1).$$
By $M{\rm j}=q^{(n-1)\ell}{\rm j}$ and $M\chi_{\{\pi\}}=v_{{\cal M}_1(\pi)}$, we obtain
$$M\chi_{\overline{{\cal M}}_n(\pi)}=(q^{-(n-1)\ell}M{\rm j}-M\chi_{\{\pi\}})
q^{{n\choose 2}}\prod_{s=1}^{n-1}(q^{\ell-n+s}-1),
$$
which implies that
$$\chi_{\overline{{\cal M}}_n(\pi)}-(q^{-(n-1)\ell}{\rm j}-\chi_{\{\pi\}})q^{{n\choose 2}}\prod_{s=1}^{n-1}(q^{\ell-n+s}-1)\in\ker(M),$$
as desired. $\qed$

A {\it partition} of the vector space $\mathbb{F}_q^{n+\ell}$ is a set ${\cal P}$ of subspaces of $\mathbb{F}_q^{n+\ell}$ such that any one-dimensional subspace is contained in exactly one element of ${\cal P}$.
If $T=\{\dim W :W\in{\cal P}\}$, the partition ${\cal P}$ is said to be a
$T$-{\it partition} of $\mathbb{F}_q^{n+\ell}$. It is well-known that an $\{n\}$-partition of $\mathbb{F}_q^{n+\ell}$ exists if and only if $n$ divides $\ell$. An $\{n\}$-partition of $\mathbb{F}_q^{n+\ell}$ is also called an $n$-{\it spread} of $\mathbb{F}_q^{n+\ell}$.

\begin{lemma}\label{lem2.8}{\rm(See \cite{Beutelspacher,Dembowski}.)}
Let $n\leq \ell$. Then the following hold:
\begin{itemize}
\item[\rm(i)]
If $n=\ell$, then $\mathbb{F}_q^{2n}$ has an $n$-spread $\{E,\pi_1,\ldots,\pi_{q^{n}}\}$, where $E=\langle e_{n+1},\ldots,e_{2n}\rangle$ and $\dim\pi_i=\cdots=\dim\pi_{q^n}=n$.

\item[\rm(ii)]
If $n<\ell$, then $\mathbb{F}_q^{n+\ell}$ has an $\{\ell,n\}$-partition $\{E,\pi_1,\ldots,\pi_{q^{\ell}}\}$, where $E=\langle e_{n+1},\ldots,e_{n+\ell}\rangle$ and $\dim\pi_i=\cdots=\dim\pi_{q^\ell}=n$.
 \end{itemize}
\end{lemma}

Let ${\cal P}=\{E,\pi_1,\ldots,\pi_{q^{\ell}}\}$ with $E=\langle e_{n+1},\ldots,e_{n+\ell}\rangle$ be an $n$-spread or an $\{\ell,n\}$-partition of $\mathbb{F}_q^{n+\ell}$ corresponding to $n=\ell$ or $n<\ell$, respectively. Then any one-dimensional subspace in the attenuated space $A_q(n+\ell,E)$ is contained in exactly one element of ${\cal P}\setminus\{E\}$.
The set ${\cal P}\setminus\{E\}$ is called an {\it attenuated $n$-spread} in $A_q(n+\ell,E)$.

Every subset of an attenuated $n$-spread in $A_q(n+\ell,E)$ is called an {\it attenuated partial $n$-spread} in $A_q(n+\ell,E)$. A {\it pair of conjugate switching sets} is a pair of disjoint attenuated partial $n$-spreads in $A_q(n+\ell,E)$ that cover the same subset of ${\cal M}_1$.

Next, we give several equivalent definitions for a Cameron-Liebler set  in ${\cal M}_n$.

\begin{thm}\label{lem2.9}
Let ${\cal L}$ be a non-empty set in ${\cal M}_n$
with $|{\cal L}|=xq^{(n-1)\ell}$, where $x$ is an integer. Then the following properties are equivalent.
\begin{itemize}
\item[\rm(i)]
$\chi_{{\cal L}}\in {\rm Im}(M^t)$.

\item[\rm(ii)]
$\chi_{{\cal L}}\in \ker(M)^{\perp}$.

\item[\rm(iii)]
For every $\pi\in{\cal M}_n$, the number of elements in ${\cal L}$ disjoint to $\pi$  is
$$(x-(\chi_{{\cal L}})_{\pi})q^{{n\choose 2}}\prod_{s=1}^{n-1}(q^{\ell-n+s}-1).$$

\item[\rm(iv)]
The vector $v=\chi_{{\cal L}}-xq^{-\ell}{\rm j}$ is a vector in $V_1$.

\item[\rm(v)]
$\chi_{{\cal L}}\in V_0\oplus V_1$.

\item[\rm(vi)]
$|{\cal L}\cap{\cal S}|=x$ for every attenuated $n$-spread ${\cal S}$ in $A_q(n+\ell,E)$.

\item[\rm(vii)]
For every pair of conjugate switching sets ${\cal R}$ and ${\cal R}'$, we have $|{\cal L}\cap{\cal R}|=|{\cal L}\cap{\cal R}'|$.
\end{itemize}
\end{thm}
\proof (i) $\Leftrightarrow$ (ii): Since ${\rm Im}(M^t)=\ker(M)^{\perp}$, the desired result follows.

(ii) $\Rightarrow$ (iii):
Let $\pi\in{\cal M}_n$. By Lemma~\ref{lem2.7}, we obtain
$$\chi_{\overline{{\cal M}}_n(\pi)}-(q^{-(n-1)\ell}{\rm j}-\chi_{\{\pi\}})q^{{n\choose 2}}\prod_{s=1}^{n-1}(q^{\ell-n+s}-1)\in\ker(M).$$
Since $\chi_{{\cal L}}\in \ker(M)^{\perp}$, we have
\begin{eqnarray*}
&&\chi_{\overline{{\cal M}}_n(\pi)}\cdot\chi_{{\cal L}}-(q^{-(n-1)\ell}{\rm j}\cdot\chi_{{\cal L}}-\chi_{\{\pi\}}\cdot\chi_{{\cal L}})q^{{n\choose 2}}\prod_{s=1}^{n-1}(q^{\ell-n+s}-1)=0\\
&\Leftrightarrow & |\overline{{\cal M}}_n(\pi)\cap{\cal L}|-(q^{-(n-1)\ell}|{\cal L}|-(\chi_{{\cal L}})_{\pi})q^{{n\choose 2}}\prod_{s=1}^{n-1}(q^{\ell-n+s}-1)=0\\
&\Leftrightarrow & |\overline{{\cal M}}_n(\pi)\cap{\cal L}|=(x-(\chi_{{\cal L}})_{\pi})q^{{n\choose 2}}\prod_{s=1}^{n-1}(q^{\ell-n+s}-1).
\end{eqnarray*}
The last equality shows that the desired result follows.

(iii) $\Rightarrow$ (iv): From (iii), we deduce that
$$K\chi_{{\cal L}}=(x{\rm j}-\chi_{{\cal L}})q^{{n\choose 2}}\prod_{s=1}^{n-1}(q^{\ell-n+s}-1)=-\lambda_1(x{\rm j}-\chi_{{\cal L}}).$$
By Lemma~\ref{lem2.6}, we have $K{\rm j}=-\lambda_1(q^\ell-1){\rm j}$, and therefore
\begin{eqnarray*}
Kv&=&K(\chi_{{\cal L}}-xq^{-\ell}{\rm j})
=-\lambda_1(x{\rm j}-\chi_{{\cal L}})+\lambda_1(q^\ell-1)xq^{-\ell}{\rm j}\\
&=&\lambda_1(\chi_{{\cal L}}+x((q^\ell-1)q^{-\ell}-1){\rm j})
=\lambda_1(\chi_{{\cal L}}-xq^{-\ell}{\rm j})
=\lambda_1v.
\end{eqnarray*}
By Lemma~\ref{lem2.4},  we obtain $v\in V_1$.

(iv) $\Rightarrow$ (v): From $V_0=\langle {\rm j}\rangle$, we deduce that the desired result follows.

(v) $\Rightarrow$ (i):  By Lemma~\ref{lem2.5}, the desired result follows.

Now we show that the property (vi) is also equivalent to the other properties.

(ii) $\Rightarrow$ (vi): Let ${\cal S}$ be an attenuated $n$-spread in $A_q(n+\ell,E)$.
Since $M\chi_{\cal S}={\rm j}$, by Lemma~\ref{lem2.1}, $\chi_{\cal S}-q^{-(n-1)\ell}{\rm j}\in\ker(M)$.  Since $\chi_{\cal L}\in\ker(M)^\bot$, we have
$$0=\chi_{\cal L}\cdot(\chi_{\cal S}-q^{-(n-1)\ell}{\rm j})=|{\cal L}\cap{\cal S}|-q^{-(n-1)\ell}|{\cal L}|,$$
which implies that $|{\cal L}\cap{\cal S}|=q^{-(n-1)\ell}|{\cal L}|=x$.

(vi) $\Rightarrow$ (iii):
Let $n_i$, for $i=1,2$, be the number
of attenuated $n$-spreads  in $A_q(n+\ell,E)$ that contain $i$ fixed pairwise disjoint subspaces in ${\cal M}_n$. Since the group $P\Gamma\!L(n+\ell,\mathbb{F}_q)_E$  acts transitively on
the couples of pairwise disjoint  subspaces in ${\cal M}_n$, this number only depends on $i$, and not on the chosen  subspaces.

For a fixed $\pi\in{\cal M}_n$, if we count the number of couples
$(\pi',{\cal S})$, where $\pi'\in{\cal M}_n$ with $\pi\cap\pi'=\{0\}$ and ${\cal S}$ is an attenuated $n$-spread  in $A_q(n+\ell,E)$ containing $\pi$ and $\pi'$, by Lemma~\ref{lem2.6}, we have
$n_1(q^\ell-1)=n_2q^{{n\choose 2}}\prod_{s=1}^{n}(q^{\ell-n+s}-1)$, which implies that $n_1/n_2=q^{{n\choose 2}}\prod_{s=1}^{n-1}(q^{\ell-n+s}-1)$.

For a fixed $\pi\in{\cal M}_n$, if we count the number of couples $(\pi',{\cal S})$, where $\pi'\in{\cal L}$ with $\pi\cap\pi'=\{0\}$ and ${\cal S}$ is an attenuated $n$-spread  in $A_q(n+\ell,E)$ containing $\pi$ and $\pi'$, then the number of subspaces in ${\cal L}$ disjoint to $\pi$ is
$$(x-(\chi_{{\cal L}})_{\pi})\frac{n_1}{n_2}=(x-(\chi_{{\cal L}})_{\pi})q^{{n\choose 2}}\prod_{s=1}^{n-1}(q^{\ell-n+s}-1).$$

Next, we show that property (vii) is equivalent with the other properties.

(ii) $\Rightarrow$ (vii):  Since  ${\cal R}$ and ${\cal R}'$ cover the same subset of ${\cal M}_1$, we have $\chi_{{\cal R}}-\chi_{{\cal R}'}\in \ker(M)$,
which implies that $\chi_{{\cal L}}\cdot(\chi_{{\cal R}}-\chi_{{\cal R}'})=\chi_{{\cal L}}\cdot\chi_{{\cal R}}-\chi_{{\cal L}}\cdot\chi_{{\cal R}'}=0$.
It follows that $|{\cal L}\cap{\cal R}|=\chi_{{\cal L}}\cdot\chi_{{\cal R}}=\chi_{{\cal L}}\cdot\chi_{{\cal R}'}=|{\cal L}\cap{\cal R}'|$.

(vii) $\Rightarrow$ (vi): For any two attenuated $n$-spreads ${\cal S}_1$ and ${\cal S}_2$ in $A_q(n+\ell,E)$,
 the the sets ${\cal S}_1\setminus{\cal S}_2$ and ${\cal S}_2\setminus{\cal S}_1$ form a pair of conjugate switching sets. So
$|{\cal L}\cap({\cal S}_1\setminus{\cal S}_2)|=|{\cal L}\cap({\cal S}_2\setminus{\cal S}_1)|$, which implies that $|{\cal L}\cap{\cal S}_1|=|{\cal L}\cap{\cal S}_2)|=c$.

Now we prove $c=x=|{\cal L}|q^{-(n-1)\ell}$.
Let $n_i$, for $i=0,1$, be the number
of attenuated $n$-spreads  in $A_q(n+\ell,E)$ that contain $i$ fixed pairwise disjoint subspaces in ${\cal M}_n$. Then this number only depends on $i$, and not on the chosen  subspaces. We count the number of couples
$(\pi,{\cal S})$, where $\pi\in{\cal M}_n$ and ${\cal S}$ is an attenuated $n$-spread  in $A_q(n+\ell,E)$ containing $\pi$, then
$n_0q^\ell=n_1q^{n\ell}$, which implies that $n_0/n_1=q^{(n-1)\ell}$.

By counting the number of couples $(\pi,{\cal S})$, where $\pi\in{\cal L}$ and ${\cal S}$ is an attenuated $n$-spread in $A_q(n+\ell,E)$ containing $\pi$, then the number of subspaces in ${\cal L}\cap{\cal S}$ equals $|{\cal L}|\frac{n_1}{n_0}=|{\cal L}|q^{-(n-1)\ell}=x$.
$\qed$

\section{Properties of Cameron-Liebler sets}
In this section, we give some examples and properties of Cameron-Liebler sets in the bilinear forms graph ${\rm Bil}_q(n,\ell)$.
We begin with a useful lemma.

\begin{lemma}\label{lem3.1}
Let ${\cal L}$ and ${\cal L}'$ be two Cameron-Liebler sets in ${\rm Bil}_q(n,\ell)$ with parameters $x$ and $x'$, respectively. Then the following hold:
\begin{itemize}
\item[\rm(i)]
$0\leq x\leq q^\ell$.

\item[\rm(ii)]
The set of all subspaces in ${\cal M}_n$ not in ${\cal L}$ is a Cameron-Liebler set in  ${\rm Bil}_q(n,\ell)$ with parameter $q^\ell-x$.

\item[\rm(iii)]
If ${\cal L}\cap{\cal L}'=\emptyset$, then ${\cal L}\cup{\cal L}'$ is a Cameron-Liebler set in  ${\rm Bil}_q(n,\ell)$ with parameter $x+x'$.

\item[\rm(iv)]
If ${\cal L}'\subseteq{\cal L}$, then ${\cal L}\setminus{\cal L}'$ is a Cameron-Liebler set in  ${\rm Bil}_q(n,\ell)$  with parameter $x-x'$.
\end{itemize}
\end{lemma}

 Next, we give some examples of Cameron-Liebler sets in  ${\rm Bil}_q(n,\ell)$.

 \begin{lemma}\label{lem3.2}
Let $\pi\in{\cal M}_n$ and $V\in{\cal M}(n+\ell-1,\ell-1;n+\ell,E)$.
Then the size of ${\cal M}_n(V)\cap\overline{{\cal M}}_n(\pi)$ is $q^{{n\choose 2}}\prod_{s=0}^{n-1}(q^{\ell-n+s}-1)$ if $\pi\in{\cal M}_n(V)$, and $q^{\ell-1+{n-1\choose 2}}\prod_{s=1}^{n-1}(q^{\ell-n+s}-1)$ if $\pi\not\in{\cal M}_n(V)$.
\end{lemma}
\proof If $\pi\in{\cal M}_n(V)$, by the transitivity of $P\Gamma\!L(n+\ell,\mathbb{F}_q)_E$ on ${\cal M}_n(V)$, we may assume that $\pi=\langle e_1,e_2,\ldots,e_{n}\rangle$ and
$V=\langle e_1,e_2,\ldots,e_{n+\ell-1}\rangle$.  Then $\pi'\in{\cal M}_n(V)\cap\overline{{\cal M}}_n(\pi)$  has a matrix representation
$$\bordermatrix{ &\hbox{\footnotesize{$n$}}\cr
&I\cr
&A\cr
&0}
\hspace{-3pt}
\begin{array}{c}
\hbox{\footnotesize{$n$}}\\
\hbox{\footnotesize{$\ell-1$}}\\
\hbox{\footnotesize{$1$}}
\end{array}\;\hbox{with rank}(A)=n.$$
By Lemma~\ref{lem2.5.2},
the size of ${\cal M}_n(V)\cap\overline{{\cal M}}_n(\pi)$ is $q^{{n\choose 2}}\prod_{s=0}^{n-1}(q^{\ell-n+s}-1)$.

If $\pi\not\in{\cal M}_n(V)$, by the transitivity of $P\Gamma\!L(n+\ell,\mathbb{F}_q)_E$ on ${\cal M}_n\setminus{\cal M}_n(V)$, we may assume that $\pi=\langle e_1+e_{n+\ell},e_2,\ldots,e_{n}\rangle$ and
$V=\langle e_1,e_2,\ldots,e_{n+\ell-1}\rangle$.  Then $\pi'\in{\cal M}_n(V)\cap\overline{{\cal M}}_n(\pi)$
 has a matrix representation
$$\bordermatrix{&\hbox{\footnotesize{$1$}} &\hbox{\footnotesize{$n-1$}}\cr
&1&0\cr
&0&I\cr
&A_1&A_2\cr
&0&0}
\hspace{-3pt}
\begin{array}{c}
\hbox{\footnotesize{$1$}}\\
\hbox{\footnotesize{$n-1$}}\\
\hbox{\footnotesize{$\ell-1$}}\\
\hbox{\footnotesize{$1$}}
\end{array}\;\hbox{with rank}(A_2)=n-1.$$
By Lemma~\ref{lem2.5.2} again, the size of ${\cal M}_n(V)\cap\overline{{\cal M}}_n(\pi)$ is
$q^{\ell-1+{n-1\choose 2}}\prod_{s=1}^{n-1}(q^{\ell-n+s}-1)$. $\qed$

Similar to Examples~3.2 and~3.3 in \cite{Blokhuis}, we obtain the following result.

\begin{lemma}\label{lem3.3}
\begin{itemize}
\item[\rm(i)]
For a fixed subspace $\tau\in{\cal M}_1$, the set ${\cal M}'_n(\tau)$ is a Cameron-Liebler set in  ${\rm Bil}_q(n,\ell)$ with parameter $1$.

\item[\rm(ii)]
For a fixed $V\in{\cal M}(n+\ell-1,\ell-1;n+\ell,E)$, the set ${\cal M}_n(V)$ is a Cameron-Liebler set in  ${\rm Bil}_q(n,\ell)$ with parameter $q^{\ell-n}$.

\item[\rm(iii)]
For two fixed subspaces $\tau\in{\cal M}_1$ and $V\in{\cal M}(n+\ell-1,\ell-1;n+\ell,E)$ with $\tau\not\subseteq V$, the set ${\cal M}'_n(\tau)\cup{\cal M}_n(V)$ is a Cameron-Liebler set in  ${\rm Bil}_q(n,\ell)$ with parameter $q^{\ell-n}+1$.

\item[\rm(iv)]
For each integer $x$ with $0\leq x\leq q^\ell$, there exists a Cameron-Liebler set in  ${\rm Bil}_q(n,\ell)$ with parameter $x$.
\end{itemize}
\end{lemma}
\proof (i). Since the characteristic vector $\chi_{{\cal M}'_n(\tau)}$ is the row of $M$ corresponding to the subspace $\tau$, by Theorem~\ref{lem2.9} (i), the set ${\cal M}'_n(\tau)$ is a Cameron-Liebler set in  ${\rm Bil}_q(n,\ell)$ with parameter $1$.

(ii). For every $\pi\in{\cal M}_n$, by Lemma~\ref{lem3.2}, the size of ${\cal M}_n(V)\cap\overline{{\cal M}}_n(\pi)$ is
$$(q^{\ell-n}-(\chi_{{\cal M}_n(V)})_{\pi})q^{{n\choose 2}}\prod_{s=1}^{n-1}(q^{\ell-n+s}-1).$$
By Theorem~\ref{lem2.9} (iii), the set ${\cal M}_n(V)$ is a Cameron-Liebler set in  ${\rm Bil}_q(n,\ell)$ with parameter $q^{\ell-n}$.

(iii).  Since $\tau\not\subseteq V$, we have ${\cal M}'_n(\tau)\cap{\cal M}_n(V)=\emptyset$, which implies that ${\cal M}'_n(\tau)\cup{\cal M}_n(V)$ is a Cameron-Liebler set in  ${\rm Bil}_q(n,\ell)$ with parameter $q^{\ell-n}+1$ by Lemma~\ref{lem3.1} (iii).

(iv). First note that a Cameron-Liebler set in ${\rm Bil}_q(n,\ell)$ with parameter 0 is the empty set.
For each $(x_1,\ldots,x_{\ell})^t\in\mathbb{F}_q^{\ell}$,
let $\tau_{(x_1,\ldots,x_{\ell})}=\langle e_1+\sum_{i=1}^{\ell}x_ie_{n+i}\rangle$.
 By (i), the set ${\cal M}'_n(\tau_{(x_1,\ldots,x_{\ell})})$ is a Cameron-Liebler set in  ${\rm Bil}_q(n,\ell)$ with parameter $1$. For any two different vectors $(x_1,\ldots,x_{\ell})^t$ and $(y_1,\ldots,y_{\ell})^t$ in $\mathbb{F}_q^{\ell}$, we have ${\cal M}'_n(\tau_{(x_1,\ldots,x_{\ell})})\cap{\cal M}'_n(\tau_{(y_1,\ldots,y_{\ell})})=\emptyset$.
 Let $S$ be a non-empty subset of $\mathbb{F}_q^{\ell}$ with size $x$.
  By Lemma~\ref{lem3.1} (iii),
the set $\bigcup_{(x_1,\ldots,x_{\ell})^t\in S}{\cal M}'_n(\tau_{(x_1,\ldots,x_{\ell})})$ is a Cameron-Liebler set in  ${\rm Bil}_q(n,\ell)$ with parameter $x$. $\qed$

Next, we give some no-trivial Cameron-Liebler sets in  ${\rm Bil}_q(n,\ell)$.
We will use the following result, the so-called Erd\H{o}s-Ko-Rado theorem for ${\rm Bil}_q(n,\ell)$.

\begin{thm}\label{lem4.1}{\rm(See \cite{Tanaka}.)}
Suppose that $\ell\geq n$ and  ${\cal F}\subseteq{\cal M}_n$ is an intersecting family. Then $|{\cal F}|\leq q^{(n-1)\ell}$, and equality holds if and only if either {\rm(a)} ${\cal F}=\{\pi\in{\cal M}_n: \tau\subseteq\pi\}$ for some $\tau\in{\cal M}_1$, or {\rm(b)} $n=\ell$ and ${\cal F}=\{\pi\in{\cal M}_n: \pi\subseteq U\}$ for some $U\in{\cal M}(2n-1,n-1;n+\ell,E)$.
\end{thm}

\begin{lemma}\label{lem3.4.0}
Let $\ell>n\geq2$. Then the following hold:
\begin{itemize}
\item[\rm(i)]
For a fixed $V\in{\cal M}(n+\ell-1,\ell-1;n+\ell,E)$, the set ${\cal M}_n(V)$ is a non-trivial Cameron-Liebler set in  ${\rm Bil}_q(n,\ell)$ with parameter $q^{\ell-n}$.

\item[\rm(ii)]
Let $1\leq y<q^{n-1}$, and $\{V_1,\ldots,V_y\}\subseteq{\cal M}(n+\ell-1,\ell-1;n+\ell,E)$ with ${\cal M}_n(V_i)\cap{\cal M}_n(V_j)=\emptyset$ for all $i\not=j$.\footnote{For each $(x_1,\ldots,x_{n})^t\in\mathbb{F}_q^{n}$, let $V_{(x_1,\ldots,x_{n})}=\langle e_1+x_1e_{n+1},\ldots,e_n+x_ne_{n+1},e_{n+2},\ldots,e_{n+\ell}\rangle.$
 Then $\{V_{(x_1,\ldots,x_{n})}: (x_1,\ldots,x_{n})^t\in\mathbb{F}_q^{n}\}\subseteq{\cal M}(n+\ell-1,\ell-1;n+\ell,E)$
and ${\cal M}_n(V_{(x_1,\ldots,x_{n})})\cap{\cal M}_n(V_{(y_1,\ldots,y_{n})})=\emptyset$
for all $(x_1,\ldots,x_{n})\not=(y_1,\ldots,y_{n})$. }
 Then the set $\bigcup_{i=1}^y{\cal M}_n(V_i)$ is a non-trivial Cameron-Liebler set in  ${\rm Bil}_q(n,\ell)$ with parameter $yq^{\ell-n}$.
\end{itemize}
\end{lemma}
\proof  (i). By Lemma~\ref{lem3.3} (ii), the set ${\cal M}_n(V)$ is a Cameron-Liebler set in  ${\rm Bil}_q(n,\ell)$ with parameter $q^{\ell-n}$. By Theorem~\ref{lem4.1}, the size of the maximum intersecting family in
${\cal M}_n(V)$ is at most $q^{(n-1)(\ell-1)}$. Since the size of ${\cal M}_n(V)$ is $q^{n(\ell-1)}$,  ${\cal M}_n(V)$ is a union of at least $q^{\ell-1}$ intersecting families,
which implies that ${\cal M}_n(V)$ is non-trivial by $q^{\ell-n}<q^{\ell-1}$.

(ii). By Lemma~\ref{lem3.1} (iii),
the set $\bigcup_{i=1}^y{\cal M}_n(V_i)$ is a Cameron-Liebler set in  ${\rm Bil}_q(n,\ell)$ with parameter $yq^{\ell-n}$. We claim that the size of the maximum intersecting family in $\bigcup_{i=1}^y{\cal M}_n(V_i)$
is at most $yq^{(n-1)(\ell-1)}$. By (i), the result is true for $y=1$.
Suppose that ${\cal F}\subseteq\bigcup_{i=1}^y{\cal M}_n(V_i)$ is a intersecting family. Then both
${\cal F}\cap\bigcup_{i=1}^{y-1}{\cal M}_n(V_i)$ and ${\cal F}\cap{\cal M}_n(V_y)$ are intersecting families.
By the induction hypothesis, we have $|{\cal F}\cap\bigcup_{i=1}^{y-1}{\cal M}_n(V_i)|\leq(y-1)q^{(n-1)(\ell-1)}$
and $|{\cal F}\cap{\cal M}_n(V_y)|\leq q^{(n-1)(\ell-1)}$, which imply that
$|{\cal F}|\leq yq^{(n-1)(\ell-1)}$. Since the size of $\bigcup_{i=1}^y{\cal M}_n(V_i)$ is $yq^{n(\ell-1)}$,  $\bigcup_{i=1}^y{\cal M}_n(V_i)$ is a union of at least $q^{\ell-1}$ intersecting families,
which implies that ${\cal M}_n(V)$ is non-trivial by $yq^{\ell-n}<q^{\ell-1}$.
$\qed$

De Boeck \cite{De Boeck3} obtained several Anzahl theorems for $n$-dimensional disjoint subspaces in the vector space $\mathbb{F}_q^{2n}$.
Now we generalize some of these results by using method in \cite{De Boeck3}.

\begin{lemma}\label{lem3.4.1}
Let $1\leq k\leq m\leq n$. Suppose that $\pi_1$ and $\pi_2$ are two disjoint $n$-dimensional subspaces in $\mathbb{F}_q^{2n}$, and
$\sigma_0$ is a $k$-dimensional subspace  in $\mathbb{F}_q^{2n}$ disjoint to both $\pi_1$ and $\pi_2$.
Then number of  $m$-dimensional subspaces in $\mathbb{F}_q^{2n}$ through $\sigma_0$ disjoint to both $\pi_1$ and $\pi_2$ equals
$$d_{km}(q,n)=q^{\frac{(m-k)(m+k-1)}{2}}{n-k\brack m-k}_q\prod_{i=1}^{m-k}(q^{n-k-i+1}-1).$$
\end{lemma}
\proof An $m$-dimensional subspace  through $\sigma_0$ is spanned by $m-k$ linearly independent vectors
not in $\sigma_0$. The first $i-1$ chosen vectors determine a $(k+i-1)$-dimensional subspace $\sigma_{i-1}$ satisfying
$\sigma_0\subseteq\sigma_{i-1}$, meeting neither $\pi_1$ nor $\pi_2$, $i=1,\ldots,m-k$. Then the $i$th vector must
be a vector not in $\langle\pi_1,\sigma_{i-1}\rangle\cup\langle\pi_2,\sigma_{i-1}\rangle$. Note that $\langle\pi_1,\sigma_{i-1}\rangle\cap\langle\pi_2,\sigma_{i-1}\rangle$ is a
$2(k+i-1)$-dimensional subspace. Hence, there are
$$q^{2n}-2q^{n+k+i-1}+q^{2(k+i-1)}=q^{2(k+i-1)}(q^{n-k-i+1}-1)^2$$
 different vectors that can be chosen as $i$th vector. So, there are
 $$\prod_{i=1}^{m-k}q^{2(k+i-1)}(q^{n-k-i+1}-1)^2$$
 different tuples of $m-k$ vectors, together with $\sigma_0$, spanning an $m$-dimensional subspace through $\sigma_0$ disjoint to both $\pi_1$ and $\pi_2$. The $m$-dimensional subspace  is defined by
$$\prod_{i=1}^{m-k}(q^m-q^{k+i-1})=\prod_{i=1}^{m-k}q^{k+i-1}(q^{m-k-i+1}-1)$$
different tuples of $m-k$ vectors. Consequently, there are
$$\prod_{i=1}^{m-k}q^{k+i-1}\frac{(q^{n-k-i+1}-1)^2}{q^{m-k-i+1}-1}=q^{\frac{(m-k)(m+k-1)}{2}}{n-k\brack m-k}_q\prod_{i=1}^{m-k}(q^{n-k-i+1}-1)$$
 such $m$-dimensional subspaces. $\qed$

 \begin{lemma}\label{lem3.4.2}
 Suppose that $1\leq k\leq m\leq n$, $\pi_1,\pi_2$ and $\pi_3$ are three pairwise disjoint $n$-dimensional subspaces in $\mathbb{F}_q^{2n}$. Let ${\cal M}_{km}(q,n)$ be the set of all $m$-dimensional subspaces  in $\mathbb{F}_q^{2n}$ disjoint to both $\pi_1$ and $\pi_2$ and meeting $\pi_3$ in a $k$-dimensional subspace, and let $x_{km}(q,n)$ be the number of pairs $(\tau,\sigma)$,
with  $\tau$ a $k$-dimensional subspace contained in $\pi_3$, $\sigma$ an $m$-dimensional subspace disjoint to both $\pi_1$ and $\pi_2$, and such that $\tau\subseteq\sigma$. Then
\begin{equation}\label{equa2}
z_{km}(q,n):=|{\cal M}_{km}(q,n)|=\sum_{i=k}^m(-1)^{i-k}{i\brack k}_qq^{{i-k\choose2}}x_{im}(q,n).
\end{equation}
Moreover, the number of $m$-dimensional subspaces in $\mathbb{F}_q^{2n}$
 disjoint to $\pi_1,\pi_2$ and $\pi_3$ equals
$$z_{0m}(q,n)=q^{{m\choose2}}{n\brack m}_q\sum_{i=0}^{m}(-1)^{i}{m\brack i}_q\prod_{j=1}^{m-i}(q^{n-i-j+1}-1),\;\hbox{where}\; \prod\emptyset=1.$$
 \end{lemma}
 \proof By counting the pairs $(\tau,\sigma)$
satisfying  $\tau\in{\cal M}_k,\sigma\in{\cal M}_{im}(q,n)$ and $\tau\subseteq\sigma\cap\pi_3$, we obtain
$$x_{km}(q,n)=\sum_{i=k}^mz_{im}(q,n){i\brack k}_q=z_{km}(q,n)+\sum_{i=k+1}^mz_{im}(q,n){i\brack k}_q.$$
Now, we prove the desired equality (\ref{equa2}) by using  induction on $t=m-k+1$. It is obvious that the equality
(\ref{equa2}) is true for $t=1$. Assume that the equality
(\ref{equa2}) is true for all $t\leq m-k$. Then
 \begin{eqnarray*}
 z_{km}(q,n)&=&x_{km}(q,n)-\sum_{i=k+1}^mz_{im}(q,n){i\brack k}_q\\
   &=&x_{km}(q,n)-\sum_{i=k+1}^m{i\brack k}_q\left(\sum_{j=i}^m(-1)^{j-i}{j\brack i}_qq^{{j-i\choose2}}x_{jm}(q,n)\right)  \\
   &=&x_{km}(q,n)-\sum_{i=k+1}^m\sum_{j=i}^m(-1)^{j-i}{i\brack k}_q{j\brack i}_qq^{{j-i\choose2}}x_{jm}(q,n)\\
   &=&x_{km}(q,n)-\sum_{j=k+1}^m(-1)^{j-k}{j\brack k}_qx_{jm}(q,n)\sum_{i=k+1}^j(-1)^{k-i}{j-k\brack i-k}_qq^{{j-i\choose2}}\\
   &=&x_{km}(q,n)+\sum_{j=k+1}^m(-1)^{j-k}{j\brack k}_qq^{{j-k\choose 2}}x_{jm}(q,n)\\
   &=&\sum_{j=k}^m(-1)^{j-k}{j\brack k}_qq^{{j-k\choose 2}}x_{jm}(q,n),
 \end{eqnarray*}
 since $\sum_{i=k+1}^j(-1)^{k-i}{j-k\brack i-k}_qq^{{j-i\choose2}}=-q^{{j-k\choose 2}}$ (cf. \cite{De Boeck3}).
 Therefore, the equality
(\ref{equa2}) is true for $t=m-k+1$.

 Since $\pi_3$ contains ${n\brack i}_q$
different $i$-dimensional subspaces, by Lemma~\ref{lem3.4.1},   $x_{im}(q,n)={n\brack i}_qd_{im}(q,n)$.  Hence, we have
\begin{eqnarray*}
z_{km}(q,n)&=&\sum_{i=k}^m(-1)^{i-k}{i\brack k}_qq^{{i-k\choose2}}x_{im}(q,n)\\
&=&\sum_{i=k}^{m}(-1)^{i-k}{i\brack k}_qq^{{i-k\choose2}}{n\brack i}_qd_{im}(q,n)\\
&=&\sum_{i=k}^{m}(-1)^{i-k}q^{{i-k\choose2}+\frac{(m-i)(m+i-1)}{2}}{i\brack k}_q{n\brack i}_q{n-i\brack m-i}_q\prod_{j=1}^{m-i}(q^{n-i-j+1}-1)\\
&=&q^{{m\choose2}}{n\brack m}_q\sum_{i=k}^{m}(-1)^{i-k}{i\brack k}_q{m\brack i}_q\prod_{j=1}^{m-i}(q^{n-i-j+1}-1)\\
\end{eqnarray*}
Note that the number of $m$-dimensional subspaces in $\mathbb{F}_q^{2n}$ disjoint to $\pi_1,\pi_2$ and $\pi_3$ is $z_{0m}(q,n)$. Consequently,
$z_{0m}(q,n)=q^{{m\choose2}}{n\brack m}_q\sum_{i=0}^{m}(-1)^{i}{m\brack i}_q\prod_{j=1}^{m-i}(q^{n-i-j+1}-1)$. $\qed$

Blokhuis, De Boeck and D'haeseleer \cite{Blokhuis}  obtained several properties of Cameron-Liebler sets of $k$-spaces in the projective space PG$(n,q)$. Next, we generalize some of these results in $A_q(n+\ell,E)$.

\begin{lemma}\label{lem3.4}
Let $\pi$ and $\pi'$ be two disjoint subspaces in ${\cal M}_n$ with $\Sigma=\langle\pi,\pi'\rangle$,
 $\tau\in{\cal M}_1(\Sigma)\setminus({\cal M}_1(\pi)\cup{\cal M}_1(\pi'))$ and $\tau'\in{\cal M}_1\setminus{\cal M}_1(\Sigma)$. For $i=0,1,\ldots,n$, let $W_i(q,n,\ell)$ be the number of subspaces $\pi''\in\overline{{\cal M}}_n(\pi)\cap\overline{{\cal M}}_n(\pi')$ with $\dim(\pi''\cap\Sigma)=i$. Then
 $$W_i(q,n,\ell)=z_{0i}(q,n)q^{n(n-i)+{n-i\choose 2}}\prod_{s=1}^{n-i}(q^{\ell-2n+i+s}-1),
$$
 where $z_{0i}(q,n)$ is given by Lemma~\ref{lem3.4.2}. Moreover,
 the size of the set $\overline{{\cal M}}_n(\pi)\cap\overline{{\cal M}}_n(\pi')$ is $W(q,n,\ell)$, the size of ${\cal M}'_n(\tau)\cap\overline{{\cal M}}_n(\pi)\cap\overline{{\cal M}}_n(\pi')$  is $W_{\Sigma}(q,n,\ell)$, and  the size of  ${\cal M}'_n(\tau')\cap\overline{{\cal M}}_n(\pi)\cap\overline{{\cal M}}_n(\pi')$ is $W_{\overline{\Sigma}}(q,n,\ell)$,
 where
 \begin{eqnarray*}
 W(q,n,\ell)&=&\sum_{i=0}^nW_i(q,n,\ell),\\
 W_{\Sigma}(q,n,\ell)&=&\frac{1}{(q^{n}-1)(q^n-2)}\sum_{i=1}^nW_i(q,n,\ell)(q^i-1),\\
 W_{\overline{\Sigma}}(q,n,\ell)&=&\frac{1}{q^n(q^{\ell-n}-1)(q^n-1)}\sum_{i=0}^{n-1}W_i(q,n,\ell)(q^n-q^i).
 \end{eqnarray*}
\end{lemma}
\proof For a given $\sigma_i\in{\cal M}_i(\Sigma)$ with $0\leq i\leq n$,  by the transitivity of $P\Gamma\!L(n+\ell,\mathbb{F}_q)_E$ on ${\cal M}_i(\Sigma)$, we may assume that $\sigma_i=\langle e_1,e_2,\ldots,e_{i}\rangle$ and
$\Sigma=\langle e_1,e_2,\ldots,e_{2n}\rangle$.  Then $\pi''\in{\cal M}_n$ with $\pi''\cap\Sigma=\sigma_i$   has a matrix representation
$$\bordermatrix{ &\hbox{\footnotesize{$i$}}&\hbox{\footnotesize{$n-i$}}\cr
&I&0\cr
&0&I\cr
&0&A_1\cr
&0&A_2}
\hspace{-3pt}
\begin{array}{c}
\hbox{\footnotesize{$i$}}\\
\hbox{\footnotesize{$n-i$}}\\
\hbox{\footnotesize{$n$}}\\
\hbox{\footnotesize{$\ell-n$}}
\end{array}\;\hbox{with rank}(A_2)=n-i.$$
By Lemma~\ref{lem2.5.2}, there are
$$
q^{n(n-i)+{n-i\choose 2}}\prod_{s=1}^{n-i}(q^{\ell-2n+i+s}-1).
$$
choices for $\pi''$ such that $\pi''\cap\Sigma=\sigma_i$.

Let ${\cal M}_{0i}(q,n,\Sigma)$ be the set of all $i$-dimensional subspaces  in $\Sigma$ disjoint to $\pi,\pi'$ and $\Sigma\cap E$. By Lemma~\ref{lem3.4.2}, $|{\cal M}_{0i}(q,n,\Sigma)|=z_{0i}(q,n)$. Since ${\cal M}_{0i}(q,n,\Sigma)\subseteq{\cal M}_i(\Sigma)$,
we obtain
$$
W_i(q,n,\ell)=z_{0i}(q,n)q^{n(n-i)+{n-i\choose 2}}\prod_{s=1}^{n-i}(q^{\ell-2n+i+s}-1)\;\hbox{for}\,i=0,1,\ldots,n,
$$
which imply  that the size of $\overline{{\cal M}}_n(\pi)\cap\overline{{\cal M}}_n(\pi')$ is
$W(q,n,\ell)=\sum_{i=0}^nW_i(q,n,\ell).$

By double counting $(\tau,\pi'')$
with $\tau\subseteq\pi''\in\overline{{\cal M}}_n(\pi)\cap\overline{{\cal M}}_n(\pi')$ and double counting
 $(\tau',\pi'')$ with $\tau'\subseteq\pi''\in\overline{{\cal M}}_n(\pi)\cap\overline{{\cal M}}_n(\pi')$,
 we obtain
 $$\left({2n\brack1}_q-3{n\brack1}_q\right)W_{\Sigma}(q,n,\ell)=\sum_{i=1}^nW_i(q,n,\ell){i\brack 1}_q\quad\hbox{and}$$
$$\left({n+\ell\brack1}_q-{\ell\brack1}_q-{2n\brack1}_q+{n\brack1}_q\right)
W_{\overline{\Sigma}}(q,n,\ell)=\sum_{i=0}^{n-1}W_i(q,n,\ell)\left({n\brack1}_q-{i\brack 1}_q\right),$$
which imply that
$$W_{\Sigma}(q,n,\ell)=\frac{1}{(q^n-1)(q^n-2)}\sum_{i=1}^nW_i(q,n,\ell)(q^i-1) \quad\hbox{and}$$
$$W_{\overline{\Sigma}}(q,n,\ell)=\frac{1}{q^n(q^{\ell-n}-1)(q^n-1)}\sum_{i=0}^{n-1}W_i(q,n,\ell)(q^n-q^i).$$
So, the proof is completed. $\qed$

For convenience, we write $W_i(q,n,\ell),W_{\Sigma}(q,n,\ell)$ and $W_{\overline{\Sigma}}(q,n,\ell)$ as
$W_i,W_{\Sigma}$ and $W_{\overline{\Sigma}}$, respectively.

\begin{lemma}\label{lem3.5}
Let ${\cal L}$ be a Cameron-Liebler set in the bilinear forms graph ${\rm Bil}_q(n,\ell)$  with parameter $x$.
Then the following hold:
\begin{itemize}
\item[\rm(i)]
For every $\pi\in{\cal L}$, there are $s_1(q,n,\ell,x)$ elements in ${\cal L}$ meeting $\pi$, where $s_1(q,n,\ell,x)=xq^{(n-1)\ell}-(x-1)q^{{n\choose 2}}\prod_{s=1}^{n-1}(q^{\ell-n+s}-1)$.

\item[\rm(ii)]
For disjoint $\pi,\pi'\in{\cal L}$ and an attenuated $n$-spread ${\cal S}_0$ in $A_q(n+n,\Sigma\cap E)$, where $\Sigma=\langle\pi,\pi'\rangle$, there are $d_2(q,n,\ell,x,{\cal S}_0)$ elements  in ${\cal L}$ disjoint to both $\pi$ and $\pi'$, and there are $s_2(q,n,\ell,x,{\cal S}_0)$
elements in ${\cal L}$ meeting both  $\pi$ and $\pi'$, where
\begin{eqnarray*}
d_2(q,n,\ell,x,{\cal S}_0)&=&(W_{\Sigma}-W_{\overline{\Sigma}})|{\cal S}_0\cap{\cal L}|-2W_{\Sigma}+xW_{\overline{\Sigma}},\\
s_2(q,n,\ell,x,{\cal S}_0)&=&xq^{(n-1)\ell}-2(x-1)q^{{n\choose 2}}\prod_{s=1}^{n-1}(q^{\ell-n+s}-1)
+d_2(q,n,\ell,x,{\cal S}_0).
\end{eqnarray*}

\item[\rm(iii)]
Let $d_2'(q,n,\ell,x)=(x-2)W_{\Sigma}$ and $s_2'(q,n,\ell,x)=xq^{(n-1)\ell}-2(x-1)q^{{n\choose 2}}\prod_{s=1}^{n-1}(q^{\ell-n+s}-1)
+d_2'(q,n,\ell,x)$. If $\ell\geq2n$, then $d_2(q,n,\ell,x,{\cal S}_0)\leq d_2'(q,n,\ell,x)$
and $s_2(q,n,\ell,x,{\cal S}_0)\leq s_2'(q,n,\ell,x)$.
\end{itemize}
\end{lemma}
\proof (i). From Theorem~\ref{lem2.9} (iii) and $|{\cal L}|=xq^{(n-1)\ell}$, we deduce that
$s_1(q,n,\ell,x)=|{\cal L}|-(x-1)q^{{n\choose 2}}\prod_{s=1}^{n-1}(q^{\ell-n+s}-1)$.

(ii).
Since $M\chi_{\{\pi\}}=v_{{\cal M}_1(\pi)}$ for each $\pi\in{\cal M}_n$,
we have $$M\chi_{{\cal S}_0}=M\sum_{\pi\in{\cal S}_0}\chi_{\{\pi\}}=\sum_{\pi\in{\cal S}_0}M\chi_{\{\pi\}}=\sum_{\pi\in{\cal S}_0}v_{{\cal M}_1(\pi)}=v_{{\cal M}_1(\Sigma)}.$$
By Lemma~\ref{lem3.4}, we obtain
\begin{eqnarray*}
&&M\chi_{\overline{{\cal M}}_n(\pi)\cap\overline{{\cal M}}_n(\pi')}\\
&=&W_{\Sigma}(v_{{\cal M}_1(\Sigma)}-v_{{\cal M}_1(\pi)}-v_{{\cal M}_1(\pi')})+W_{\overline{\Sigma}}({\rm j}-v_{{\cal M}_1(\Sigma)})\\
&=&W_{\Sigma}(M\chi_{{\cal S}_0}-M\chi_{\{\pi\}}-M\chi_{\{\pi'\}})+W_{\overline{\Sigma}}(q^{-(n-1)\ell}M{\rm j}-
M\chi_{{\cal S}_0})\\
&\Leftrightarrow&\chi_{\overline{{\cal M}}_n(\pi)\cap\overline{{\cal M}}_n(\pi')}-W_{\Sigma}(\chi_{{\cal S}_0}-\chi_{\{\pi\}}-\chi_{\{\pi'\}})-W_{\overline{\Sigma}}(q^{-(n-1)\ell}{\rm j}-
\chi_{{\cal S}_0})\in\ker(M).
\end{eqnarray*}
Since $\chi_{{\cal L}}\in\ker(M)^\bot$, we have
\begin{eqnarray*}
&&\chi_{\overline{{\cal M}}_n(\pi)\cap\overline{{\cal M}}_n(\pi')}\cdot\chi_{{\cal L}}
=W_{\Sigma}(\chi_{{\cal S}_0}\cdot\chi_{{\cal L}}-(\chi_{{\cal L}})_{\pi}-(\chi_{{\cal L}})_{\pi'})+W_{\overline{\Sigma}}(x-\chi_{{\cal S}_0}\cdot\chi_{{\cal L}})\\
&\Leftrightarrow& |\overline{{\cal M}}_n(\pi)\cap\overline{{\cal M}}_n(\pi')\cap{\cal L}|
=W_{\Sigma}(|{\cal S}_0\cap{\cal L}|-2)+W_{\overline{\Sigma}}(x-|{\cal S}_0\cap{\cal L}|)\\
&\Leftrightarrow& |\overline{{\cal M}}_n(\pi)\cap\overline{{\cal M}}_n(\pi')\cap{\cal L}|
=(W_{\Sigma}-W_{\overline{\Sigma}})|{\cal S}_0\cap{\cal L}|-2W_{\Sigma}+xW_{\overline{\Sigma}},
\end{eqnarray*}
which gives the formula for $d_2(q,n,\ell,x,{\cal S}_0)$.
The formula for $s_2(q,n,\ell,x,{\cal S}_0)$ follows from the
inclusion-exclusion principle.

(iii). Suppose that $\Sigma\in{\cal M}(2n,n;n+\ell,E)$ and ${\cal S}_0$ is an attenuated $n$-spread in  $A_q(n+n,\Sigma\cap E)$ such that $|{\cal S}_0\cap{\cal L}|>x$.
By Theorem~\ref{lem2.9} (i), we obtain that $\chi_{{\cal L}}$ can be written
as $\sum_{\tau\in{\cal M}_1}x_{\tau}\chi_{{\cal M}'_n(\tau)}$ for some $x_{\tau}\in\mathbb{R}$.
 From $\chi_{{\cal L}}\cdot {\rm j}=xq^{(n-1)\ell}$, we deduce that
$\sum_{\tau\in{\cal M}_1}x_{\tau}=x$. Since $|{\cal S}_0\cap{\cal L}|=\chi_{{\cal L}}\cdot\chi_{{\cal S}_0}=\sum_{\tau\in{\cal M}_1(\Sigma)}x_{\tau}>x$, we have
$$\sum_{\tau\in{\cal M}_1\setminus{\cal M}_1(\Sigma)}x_{\tau}=\sum_{\tau\in{\cal M}_1}x_{\tau}-\sum_{\tau\in{\cal M}_1(\Sigma)}x_{\tau}<0.$$

For every $\pi\in{\cal M}_n$, we obtain that $\chi_{{\cal L}}\cdot\chi_{\{\pi\}}=\sum_{\tau\in{\cal M}_1(\pi)}x_{\tau}$
equals 1 if $\pi\in{\cal L}$ and 0 otherwise.
For each $\tau'\in{\cal M}_1\setminus{\cal M}_1(\Sigma)$, by $\ell\geq2n$ and the proof of Lemma~\ref{lem3.4},
 there are  $$\frac{W_0(q,n,\ell)}{q^n(q^{\ell-n}-1)}=q^{3{n\choose 2}}\prod_{s=1}^{n-1}(q^{\ell-2n+s}-1)$$ choices for  $\pi'\in{\cal M}_n$ such that $\pi'\cap\Sigma=\{0\}$ and $\tau'\subseteq\pi'$. Therefore, we obtain
 $$\sum_{\tau'\in{\cal M}_1\setminus{\cal M}_1(\Sigma)}x_{\tau'}=\frac{1}{q^{3{n\choose 2}}\prod_{s=1}^{n-1}(q^{\ell-2n+s}-1)}\sum_{\pi'\in{\cal M}_n {\rm with} \pi'\cap\Sigma=\{0\}}\sum_{\tau'\in{\cal M}_1(\pi')}x_{\tau'}\geq0,$$
a contradiction.

Therefore, we have $|{\cal S}_0\cap{\cal L}|\leq x$. Since this is true for every attenuated $n$-spread ${\cal S}_0$ in  $A_q(n+n,\Sigma\cap E)$ for every $\Sigma\in{\cal M}(2n,n;n+\ell,E)$, we have
$d_2(q,n,\ell,x,{\cal S}_0)\leq(x-2)W_{\Sigma}=d'_2(q,n,\ell,x)$ by $W_{\Sigma}\geq W_{\overline{\Sigma}}$.
So, $s_2(q,n,\ell,x,{\cal S}_0)\leq s'_2(q,n,\ell,x)$. $\qed$

For convenience, we write $s_1(q,n,\ell,x),s_2(q,n,\ell,x,{\cal S}_0)$ and $s'_2(q,n,\ell,x)$ as
$s_1,s_2$ and $s_2'$, respectively.

\begin{lemma}\label{lem3.6}
Let $c,n,\ell$ be nonnegative integers with $\ell\geq2n$ and
$$(c+1)s_1-{c+1\choose 2}s_2'\geq xq^{(n-1)\ell}.$$
Then no Cameron-Liebler set  in the bilinear forms graph ${\rm Bil}_q(n,\ell)$ with parameter $x$ contains $c+1$ mutually disjoint subspaces.
\end{lemma}
\proof Suppose that ${\cal L}$ is a Cameron-Liebler set in ${\rm Bil}_q(n,\ell)$ with parameter $x$ that
contains $c+1$ mutually disjoint subspaces $\pi_0,\pi_1,\ldots,\pi_c$. By Lemma~\ref{lem3.5}, $\pi_i$ meets at least
$s_1-is_2$ elements of ${\cal L}$ that are disjoint to $\pi_0,\ldots,\pi_{i-1}$. Hence
$$xq^{(n-1)\ell}=|{\cal L}|\geq(c+1)s_1-s_2\sum_{i=1}^{c}i=(c+1)s_1-{c+1\choose 2}s_2\geq(c+1)s_1-{c+1\choose 2}s'_2,$$ which contradicts the assumption. $\qed$

\section{Classification results}
In this section, we will list some classification results for Cameron-Liebler sets in the bilinear forms graph ${\rm Bil}_q(n,\ell)$.

Firstly, we give a classification result for Cameron-Liebler sets in ${\rm Bil}_q(n,\ell)$ with parameter $x=1$.

\begin{thm}\label{lem4.2} Let ${\cal L}$ be a Cameron-Liebler set in ${\rm Bil}_q(n,\ell)$ with parameter $x=1$.
Then ${\cal L}$ is trivial.
\end{thm}
\proof By Theorem~\ref{lem2.9} (iii), ${\cal L}$ is an intersecting family, which implies that the theorem follows by Theorem~\ref{lem4.1}. $\qed$

Let ${\cal L}$ be a Cameron-Liebler set in ${\rm Bil}_q(n,\ell)$ with parameter $x$. By Lemma~\ref{lem3.1}, we have
$0\leq x\leq q^\ell$. Suppose $\ell\geq2n\geq4$. Next, we will prove that there are no non-trivial Cameron-Liebler sets  in ${\rm Bil}_q(n,\ell)$ with parameter $2\leq x\leq (q-1)^{\frac{n}{2}}q^{\frac{\ell-2n+1}{2}}$ by using method in \cite{Blokhuis}.
We will begin with the following result, the so-called Hilton-Milner theorem for ${\rm Bil}_q(n,\ell)$.

\begin{thm}\label{lem4.3}{\rm(See \cite{Gong}.)}
Suppose that $\ell\geq n+1\geq3$ with $(q,\ell)\not=(2,n+1)$, and  ${\cal F}\subseteq{\cal M}_n$ is an intersecting family with $\bigcap_{\pi'\in{\cal F}}\pi'=\{0\}$. Then the following hold:
\begin{itemize}
\item[\rm(i)]
If $n\not=3$, then $|{\cal F}|\leq q^{(n-1)\ell}-q^{{n\choose 2}}\prod_{s=1}^{n-1}(q^{\ell-n+s}-1)+q^{n-1}(q-1)$, and equality holds if and only if ${\cal F}=\{\pi'\in{\cal M}_n: \tau\subseteq\pi',\dim(\pi\cap\pi')\geq1\}\cup\{\pi'\in{\cal M}_n:\pi'\subseteq\langle\tau,\pi\rangle\}$ for some $\tau\in{\cal M}_1$ and $\pi\in{\cal M}_n$ with $\tau\not\subseteq\pi$.

\item[\rm(ii)]
If $n=3$, then $|{\cal F}|\leq q^\ell(q^2+q+1)-q(q+1)$, and equality holds if and only if
${\cal F}=\{\pi'\in{\cal M}_3: \dim(\pi\cap\pi')\geq2\}$ for some $\pi\in{\cal M}_n$.
    \end{itemize}
\end{thm}

\begin{lemma}\label{lem4.4}
Let $\ell\geq n\geq2$. Then
$$q^{(n-1)\ell}>q^{{n\choose 2}}\prod_{s=1}^{n-1}(q^{\ell-n+s}-1)>W_{\Sigma}.$$
\end{lemma}
\proof Let $\pi$ and $\pi'$ be two disjoint subspaces in ${\cal M}_n$ with $\Sigma=\langle\pi,\pi'\rangle$ and
 $\tau\in{\cal M}_1(\Sigma)\setminus({\cal M}_1(\pi)\cup{\cal M}_1(\pi'))$, since ${\cal M}'_n(\tau)\supseteq{\cal M}'_n(\tau)\cap\overline{{\cal M}}_n(\pi)\supseteq{\cal M}'_n(\tau)\cap\overline{{\cal M}}_n(\pi)\cap\overline{{\cal M}}_n(\pi')$, the desired result follows by Lemmas~\ref{lem2.6} and~\ref{lem3.4}. $\qed$

 From now on, we denote $\Delta(q,n,\ell)=q^{{n\choose 2}}\prod_{s=1}^{n-1}(q^{\ell-n+s}-1)$ and $C(q,n,\ell)=q^{(n-1)\ell}-\Delta(q,n,\ell)$.
 By Lemma~\ref{lem3.5},  we may write
 \begin{eqnarray}
 s_1(q,n,\ell,x)&=&xC(q,n,\ell)+\Delta(q,n,\ell)\label{shi-b}\quad\hbox{and}\quad\\
 s'_2(q,n,\ell,x)&=&xC(q,n,\ell)-(x-2)\Delta(q,n,\ell)+(x-2)W_{\Sigma}.\label{shi-c}
 \end{eqnarray}

 \begin{lemma}\label{lem4.4.1}
Let $\ell\geq 2n\geq4$. If $x\leq(q-1)^{\frac{n}{2}}q^{\frac{\ell-2n+1}{2}}$, then
$$\Delta(q,n,\ell)>x^2C(q,n,\ell).$$
\end{lemma}
\proof For given $\tau\in{\cal M}_1$ and $\pi\in{\cal M}_n$ with $\tau\not\subseteq\pi$, by Lemma~\ref{lem2.6}, the size of $\overline{{\cal M}}_n(\pi)\cap{\cal M}'_n(\tau)$ is $\Delta(q,n,\ell)$. Therefore, the number of $\pi'\in{\cal M}_n$ satisfying $\tau\subseteq\pi'$ and $\pi\cap\pi'\not=\{0\}$ equals $C(q,n,\ell)$.
Since this number is smaller than the product of the number of one-dimensional subspaces $\tau'\subseteq\pi'$ and the number of $n$-dimensional subspaces $\langle\tau,\tau'\rangle\subseteq\pi'\in {\cal M}_n$. By Lemma~\ref{lem2.1}, this implies that
\begin{equation}\label{shi-z}
C(q,n,\ell)\leq{n\brack 1}_qq^{\ell(n-2)}<\frac{q^{n+\ell(n-2)}}{q-1}.
\end{equation}
From $x\leq(q-1)^{\frac{n}{2}}q^{\frac{\ell-2n+1}{2}}$ and (\ref{shi-z}), we deduce that
\begin{eqnarray*}
x^2C(q,n,\ell)
&\leq&(q-1)^nq^{\ell-2n+1}\frac{q^{n+\ell(n-2)}}{q-1}\\
&=&(q-1)^{n-1}q^{(\ell-1)(n-1)}\\
&=&q^{{n\choose 2}}\prod_{s=1}^{n-1}q^{\ell-n+s-1}(q-1)\\
&<&q^{{n\choose 2}}\prod_{s=1}^{n-1}(q^{\ell-n+s}-1),
\end{eqnarray*}
which implies that the desired result follows. $\qed$

\begin{lemma}\label{lem4.4.2}
Let $\ell\geq 2n\geq4$. Then
$$W_{\Sigma}\leq \Delta(q,n,\ell)-C(q,n,\ell).$$
\end{lemma}
\proof We divide this proof into two steps.

{\bf Step}~1:  We prove $W_{\Sigma}\leq\left(1-\frac{q^{n}-2q+1}{q(q-1)(q^{\ell-1}-1)}\right)\Delta(q,n,\ell)$.

For given  $\pi,\pi'\in{\cal M}_n$ and $\tau\in{\cal M}_1(\Sigma)\setminus({\cal M}_1(\pi)\cup{\cal M}_1(\pi'))$ with $\pi\cap\pi'=\{0\}$ and $\Sigma=\langle\pi,\pi'\rangle$, by Lemmas~\ref{lem2.6} and~\ref{lem3.4}, the number of $\pi''\in{\cal M}_n$ satisfying $\tau\subseteq\pi'',\pi\cap\pi''=\{0\}$ and $\pi\cap\pi'\not=\{0\}$ equals $q^{{n\choose 2}}\prod_{s=1}^{n-1}(q^{\ell-n+s}-1)-W_{\Sigma}$. Let $\Omega(q,n,\ell)$ be the number of $\pi''\in{\cal M}_n$ satisfying $\tau\subseteq\pi'',\pi\cap\pi''=\{0\}$ and $\dim(\pi\cap\pi')=1$. Then $\Delta(q,n,\ell)-W_{\Sigma}\geq\Omega(q,n,\ell).$

We claim that \begin{equation}\label{shi-a}\Omega(q,n,\ell)\geq\frac{q^{n}-2q+1}{q(q-1)(q^{\ell-1}-1)}\Delta(q,n,\ell).\end{equation}
By the transitivity of $P\Gamma\!L(n+\ell,\mathbb{F}_q)_E$ on ${\cal M}_n$, we may assume that $\pi=\langle e_1,e_2,\ldots,e_{n}\rangle$ and $\pi'=\langle e_1+e_{n+1},e_2+e_{n+2},\ldots,e_{n}+e_{2n}\rangle$.  Since $\tau\in{\cal M}_1(\Sigma)\setminus({\cal M}_1(\pi)\cup{\cal M}_1(\pi'))$, we may assume that $\tau=\langle\sum_{i=1}^n(x_ie_i+y_ie_{n+i})\rangle$ with $(x_1,\ldots,x_n)\not=0$ and $(y_1,\ldots,y_n)\not=0,(x_1,\ldots,x_n)$, where
$(x_1,\ldots,x_n)^t,(y_1,\ldots,y_n)^t\in\mathbb{F}_q^n$.
By the transitivity of the stabilizer of both $\pi$ and $\pi'$ in $P\Gamma\!L(n+\ell,\mathbb{F}_q)_E$,
we may assume that $x_1=1$ and $x_2=\cdots=x_n=0$. Then there are the following two cases to be considered.

{\it Case}~1: $y_2=\cdots=y_n=0$. Without loss of generality, we may assume that $y_1=\xi$ for some $\xi\in\mathbb{F}_q\setminus\{0,1\}$ with $q\geq3$.
So $\tau=\langle e_1+\xi e_{n+1}\rangle$. Note that there are $q{n-1\brack 1}_q$ choices for $\tau'\subseteq\pi'$ with $\langle\tau,\tau'\rangle\in{\cal M}_2$ and $\langle\tau,\tau'\rangle\cap\pi=\{0\}$. For any given $\tau'$, without loss of generality, we may assume that
$\tau'=\langle e_2+e_{n+2}\rangle$, there exist $(n-2)$-dimensional subspaces
$$\sigma=\bordermatrix{ &\hbox{\footnotesize{$n-2$}}\cr
&0\cr
&I\cr
&A\cr
&B}
\hspace{-3pt}
\begin{array}{c}
\hbox{\footnotesize{$2$}}\\
\hbox{\footnotesize{$n-2$}}\\
\hbox{\footnotesize{$2$}}\\
\hbox{\footnotesize{$\ell-2$}}
\end{array} \quad\hbox{with rank}\,B=n-2,$$
such that $\langle\tau,\tau',\sigma\rangle\in{\cal M}_n,\langle\tau,\tau',\sigma\rangle\cap\pi=\{0\}$ and $\langle\tau,\tau',\sigma\rangle\cap\pi'=\tau'$.
 By Lemma~\ref{lem2.5.2}, there are
 $$q^{\frac{(n+1)(n-2)} {2}}\prod_{s=1}^{n-2}(q^{\ell-n+s}-1)$$ choices for $\sigma$.
 Therefore, we obtain
 $$\Omega(q,n,\ell)\geq \frac{(q^{n}-q)q^{\frac{(n+1)(n-2)} {2}}\prod_{s=1}^{n-2}(q^{\ell-n+s}-1)}{q-1}\quad\hbox{with}\quad q\geq3.$$

 {\it Case}~2: $(y_2,\ldots,y_n)\not=0$. Without loss of generality, we may assume that $y_2=1$ and $y_1=y_3=\cdots=y_n=0$.
So $\tau=\langle e_1+e_{n+2}\rangle$. Note that there are $q^2{n-2\brack 1}_q+q-1$ choices for $\tau'\subseteq\pi'$ with $\langle\tau,\tau'\rangle\in{\cal M}_2$ and $\langle\tau,\tau'\rangle\cap\pi=\{0\}$. Similar to the proof of Case~1, for any given $\tau'$, there exist
$$q^{\frac{(n+1)(n-2)} {2}}\prod_{s=1}^{n-2}(q^{\ell-n+s}-1)$$ many $(n-2)$-dimensional subspaces $\sigma$, such that
$\langle\tau,\tau',\sigma\rangle\in{\cal M}_n,\langle\tau,\tau',\sigma\rangle\cap\pi=\{0\}$ and $\langle\tau,\tau',\sigma\rangle\cap\pi'=\tau'$.
Therefore, we obtain
$$\Omega(q,n,\ell)\geq \frac{(q^{n}-2q+1)q^{\frac{(n+1)(n-2)} {2}}\prod_{s=1}^{n-2}(q^{\ell-n+s}-1)}{q-1}.$$

Since
$$\frac{\Omega(q,n,\ell)}{\Delta(q,n,\ell)}\geq\frac{(q^{n}-2q+1)q^{\frac{(n+1)(n-2)} {2}}\prod_{s=1}^{n-2}(q^{\ell-n+s}-1)}{(q-1)\Delta(q,n,\ell)}
=\frac{q^{n}-2q+1}{q(q-1)(q^{\ell-1}-1)},$$
the desired inequality (\ref{shi-a}) follows. Since $\Delta(q,n,\ell)-W_{\Sigma}\geq \Omega(q,n,\ell)$,
we complete the proof of Step~1.

{\bf Step}~2: We prove $W_{\Sigma}\leq\Delta(q,n,\ell)-C(q,n,\ell)$.

Since $W_{\Sigma}\leq\left(1-\frac{q^{n}-2q+1}{q(q-1)(q^{\ell-1}-1)}\right)\Delta(q,n,\ell)$, we only need to prove
$$\frac{q^{n}-2q+1}{q(q-1)(q^{\ell-1}-1)}\Delta(q,n,\ell)\geq C(q,n,\ell).$$ Note that
\begin{eqnarray*}
\frac{q^{n}-2q+1}{q(q-1)(q^{\ell-1}-1)}\Delta(q,n,\ell)&=&\frac{q^{n}-2q+1}{q(q-1)(q^{\ell-1}-1)}q^{{n\choose 2}}\prod_{s=1}^{n-1}(q^{\ell-n+s}-1)\\
&\geq&\frac{q^{n}-2q+1}{q(q-1)}q^{{n\choose 2}}\prod_{s=1}^{n-2}q^{\ell-n+s-1}\\
&=&\frac{q^{n}-2q+1}{q-1}q^{\ell(n-2)+2}\geq\frac{q^{n+\ell(n-2)}}{q-1}\\
&>&C(q,n,\ell),
\end{eqnarray*}
where the final inequality is given by (\ref{shi-z}).
We complete the proof of Step~2.
 $\qed$

\begin{lemma}\label{lem4.5}
 Let $\ell\geq2n\geq4$ and ${\cal L}$ be a Cameron-Liebler set in ${\rm Bil}_q(n,\ell)$ with parameter
$2\leq x\leq (q-1)^{\frac{n}{2}}q^{\frac{\ell-2n+1}{2}}$.  Then ${\cal L}$ cannot contain $\lfloor\frac{3x}{2}\rfloor$ mutually disjoint subspaces in ${\cal M}_n$.
\end{lemma}
\proof
Suppose $2\leq x\leq (q-1)^{\frac{n}{2}}q^{\frac{\ell-2n+1}{2}}$. We claim that
$$\left\lfloor\frac{3x}{2}\right\rfloor s_1-{\lfloor\frac{3x}{2}\rfloor\choose 2}s_2'>xq^{(n-1)\ell}.$$
For convenience, we write $\Delta(q,n,\ell)$ and $C(q,n,\ell)$ as
$\Delta$ and $C$, respectively. By Lemma~\ref{lem4.4.2}, (\ref{shi-b}) and (\ref{shi-c}),
 the following inequality is sufficient:
\begin{eqnarray*}
&&\left\lfloor\frac{3x}{2}\right\rfloor(xC+\Delta)-x(C+\Delta)\\
&&\quad-\frac{1}{2}\left\lfloor\frac{3x}{2}\right\rfloor\left(\left\lfloor\frac{3x}{2}\right\rfloor-1\right)(xC-(x-2)\Delta+(x-2)(\Delta-C))>0\\
&\Leftrightarrow&\Delta\left(\left\lfloor\frac{3x}{2}\right\rfloor-x\right)>C\left(x-\left\lfloor\frac{3x}{2}\right\rfloor x
+\left\lfloor\frac{3x}{2}\right\rfloor\left(\left\lfloor\frac{3x}{2}\right\rfloor-1\right)\right).
\end{eqnarray*}
From Lemma~\ref{lem4.4.1}, we deduce that $\Delta>x^2C$, which implies that the following inequality is sufficient:
\begin{equation}\label{shi-d}
x^2\left(\left\lfloor\frac{3x}{2}\right\rfloor-x\right)>x-\left\lfloor\frac{3x}{2}\right\rfloor x
+\left\lfloor\frac{3x}{2}\right\rfloor\left(\left\lfloor\frac{3x}{2}\right\rfloor-1\right).
\end{equation}
If $x$ is even, then $\lfloor\frac{3x}{2}\rfloor=\frac{3x}{2}$, which implies that the inequality (\ref{shi-d}) is equivalent with
$$2x^3-3x^2+2x>0.$$
Since $x\geq 2$, the inequality $2x^2-3x+2=2\left(x-\frac{3}{4}\right)^2+\frac{7}{8}>0$ is sufficient.
If $x$ is odd, then $\lfloor\frac{3x}{2}\rfloor=\frac{3x-1}{2}$, which implies that  the inequality (\ref{shi-d}) is equivalent with
$$f(x)=2x^3-5x^2+6x-3>0.$$
Since $x\geq 2$ and $f(2)=5$,  the inequality $f'(x)=6x^2-10x+6=6\left(x-\frac{5}{6}\right)^2+\frac{11}{6}>0$ is sufficient.
By Lemma~\ref{lem3.6},  ${\cal L}$ cannot contain $\lfloor\frac{3x}{2}\rfloor$ mutually disjoint subspaces in ${\cal M}_n$. $\qed$

\begin{lemma}\label{lem4.6}
Let $\tau\in{\cal M}_1$ and $\pi\in{\cal M}_n$ with $\tau\not\subseteq\pi$. Then the number of subspaces $\pi'\in{\cal M}_n$ satisfying $\tau\subseteq\pi'$ and $\dim(\pi\cap\pi')=1$ equals $q^{{n-1\choose 2}+1}{n-1\brack1}_q\prod_{s=2}^{n-1}(q^{\ell-n+s}-1)$.
\end{lemma}
\proof
By the transitivity of the group $P\Gamma\!L(n+\ell,\mathbb{F}_q)_E$ on $\overline{{\cal M}}_n(\pi)$, we may assume that $\tau=\langle e_1+e_{n+1}\rangle$ and $\pi=\langle e_1,\ldots,e_n\rangle$. Note that there are ${n\brack1}_1-1$ subspaces $\tau'\in{\cal M}_1(\pi)$ such that $\langle\tau,\tau'\rangle\in{\cal M}_2$.
For a given $\tau'\in{\cal M}_1(\pi)\setminus\{\langle e_1\rangle\}$, the number of $\pi'\in{\cal M}_n$ satisfying $\tau\subseteq\pi'$ and $\pi\cap\pi'=\tau'$ does not depend on the particular choice of $\tau'$. Without loss of generality, assume that
$\tau'=\langle e_2\rangle$.
Then $\pi'\in{\cal M}_n$ satisfying $\tau\subseteq\pi'$ and $\pi\cap\pi'=\tau'$  has a matrix representation
$$\bordermatrix{ &\hbox{\footnotesize{$1$}} &\hbox{\footnotesize{$1$}}&\hbox{\footnotesize{$n-2$}}\cr
&1&0&0\cr
&0&1&0\cr
&0&0&I\cr
&1&0&A_{1}\cr
&0&0&A_2}
\hspace{-3pt}
\begin{array}{c}
\hbox{\footnotesize{$1$}}\\
\hbox{\footnotesize{$1$}}\\
\hbox{\footnotesize{$n-2$}}\\
\hbox{\footnotesize{$1$}}\\
\hbox{\footnotesize{$\ell-1$}}
\end{array}\;\hbox{with rank}(A_2)=n-2.$$
By Lemma~\ref{lem2.5.2}, the number of $\pi'\in{\cal M}_n$ satisfying $\tau\subseteq\pi'$ and $\pi\cap\pi'=\tau'$  is equal to
$$q^{n-2}q^{{n-2\choose 2}}\prod_{s=1}^{n-2}(q^{\ell-n+s+1}-1)=q^{{n-1\choose 2}}\prod_{s=2}^{n-1}(q^{\ell-n+s}-1).$$
Therefore,
the number of subspaces $\pi'\in{\cal M}_n$ satisfying $\tau\subseteq\pi'$ and $\dim(\pi\cap\pi')=1$ equals $q^{{n-1\choose 2}+1}{n-1\brack1}_q\prod_{s=2}^{n-1}(q^{\ell-n+s}-1)$. $\qed$

\begin{lemma}\label{lem4.7}
Let $\ell\geq2n\geq4$ and $2\leq x\leq (q-1)^{\frac{n}{2}}q^{\frac{\ell-2n+1}{2}}$. Then the following hold:
\begin{itemize}
\item[\rm(i)]
If $n\not=3$, then
$$
\frac{x-1}{\lfloor\frac{3x}{2}\rfloor-2}\Delta(q,n,\ell)-\left(\left\lfloor\frac{3x}{2}\right\rfloor-3\right)s'_2>C(q,n,\ell)+q^{n-1}(q-1).
$$

\item[\rm(ii)]
If $n=3$, then
 $$\frac{x-1}{\lfloor\frac{3x}{2}\rfloor-2}\Delta(q,3,\ell)-\left(\left\lfloor\frac{3x}{2}\right\rfloor-3\right)s'_2
>q^\ell(q^2+q+1)-q(q+1).$$
\end{itemize}
\end{lemma}
\proof For $n\geq 2$, we first prove the following inequality:
\begin{equation}\label{equa4}
\frac{x-1}{\lfloor\frac{3x}{2}\rfloor-2}\Delta(q,n,\ell)-\left(\left\lfloor\frac{3x}{2}\right\rfloor-3\right)s'_2
>xC(q,n,\ell).
\end{equation}
By (\ref{shi-b}) and (\ref{shi-c}), the inequality (\ref{equa4}) is equivalent with
$$\frac{x-1}{\lfloor\frac{3x}{2}\rfloor-2}\Delta(q,n,\ell)-\left(\left\lfloor\frac{3x}{2}\right\rfloor-3\right)(xC(q,n,\ell)-(x-2)\Delta(q,n,\ell)+(x-2)W_{\Sigma})>xC(q,n,\ell).$$
By Lemmas~\ref{lem4.4.1} and~\ref{lem4.4.2},
 the following inequality is sufficient:
\begin{equation}\label{shi-h}\frac{x^2(x-1)}{\lfloor\frac{3x}{2}\rfloor-2}>2\left(\left\lfloor\frac{3x}{2}\right\rfloor-3\right)+x.
\end{equation}
If $x$ is even, then the inequality (\ref{shi-h}) is equivalent with
$$g(x)=x^3-7x^2+17x-12>0.$$
Since $x\geq 2$ and $g(2)=2$,  the inequality $g'(x)=3x^2-14x+17=3\left(x-\frac{7}{3}\right)^2+\frac{2}{3}>0$ is sufficient.
If $x$ is odd, then the inequality (\ref{shi-h}) is equivalent with
$$h(x)=2x^3-14x^2+41x-35>0.$$
Since $x\geq 2$ and $h(2)=7$,  the inequality $h'(x)=6x^2-28x+41=6\left(x-\frac{7}{3}\right)^2+\frac{25}{3}>0$ is sufficient.
 Therefore, the desired inequality (\ref{equa4}) follows.

(i).  Firstly, we  prove
$(x-1)C(q,n,\ell)>q^{n-1}(q-1)$ for $n\geq 4$.
For given $\tau\in{\cal M}_1$ and $\pi\in{\cal M}_n$ with $\tau\not\subseteq\pi$, by Lemma~\ref{lem2.6}, the number of $\pi'\in{\cal M}_n$ satisfying $\tau\subseteq\pi'$ and $\pi\cap\pi'\not=\{0\}$ equals $C(q,n,\ell)$. We know that this number is larger than the number of $\pi'\in{\cal M}_n$ satisfying $\tau\subseteq\pi'$ and $\dim(\pi\cap\pi')=1$, which equals $q^{{n-1\choose 2}+1}{n-1\brack1}_q\prod_{s=2}^{n-1}(q^{\ell-n+s}-1)$ by Lemma~\ref{lem4.6}.
Since $x\geq2$ and $q^{{n-1\choose 2}+1}{n-1\brack1}_q>q^{{n\choose2}}>q^{n-1}(q-1)$ for $n\geq 4$, we find that
$$(x-1)q^{{n-1\choose 2}+1}{n-1\brack1}_q\prod_{s=2}^{n-1}(q^{\ell-n+s}-1)>q^{n-1}(q-1)$$
is sufficient.

Next, we prove
$\frac{x-1}{\lfloor\frac{3x}{2}\rfloor-2}\Delta(q,2,\ell)-\left(\left\lfloor\frac{3x}{2}\right\rfloor-3\right)s'_2>C(q,2,\ell)+q(q-1)$ for $n=2$.
Note that $\Delta(q,2,\ell)=q(q^{\ell-1}-1)$ and $C(q,2,\ell)=q$. By (\ref{shi-c}) and Lemma~\ref{lem4.4.2}, we have
\begin{eqnarray*}
s_2'&=& xq-(x-2)q(q^{\ell-1}-1)+(x-2)W_{\Sigma}\\
&\leq& xq-(x-2)q(q^{\ell-1}-1)+(x-2)q(q^{\ell-1}-2)\\
&=&2q.
\end{eqnarray*}
Therefore, the following inequality is sufficient:
\begin{equation}\label{shi-g}
\frac{x-1}{\lfloor\frac{3x}{2}\rfloor-2}(q^{\ell-1}-1)-2\left(\left\lfloor\frac{3x}{2}\right\rfloor-3\right)>q.
\end{equation}

If $x$ is even, then the inequality (\ref{shi-g}) is equivalent with
$$g(x,q)=9x^2-(2q^{\ell-1}-3q+28)x+2q^{\ell-1}-4q+22<0.$$
Let
$$
\alpha=\frac{2q^{\ell-1}-3q+28+\sqrt{\Theta}}{18} \quad\hbox{with}\quad\Theta=(2q^{\ell-1}-3q+28)^2-36(2q^{\ell-1}-4q+22).$$
Then $g(\alpha,q)=0$. Since $g(2,q)=-2q^{\ell-1}+2q+2<0$, the  inequality $(q-1)^{\frac{n}{2}}q^{\frac{\ell-2n+1}{2}}<\alpha$ is sufficient.
If $q=2$, then $\sqrt{\Theta}=\sqrt{2^{2\ell}+8\cdot2^{\ell}-20}> 2^{\ell}-2$, which implies that
$\frac{1}{9}(2^{\ell}+10)<\alpha$, and therefore the following inequality is sufficient:
\begin{equation}\label{even-1}
2^{\frac{\ell-2n+1}{2}}\leq \frac{1}{9}(2^{\ell}+10).
\end{equation}
If $\ell=2n=4$, the inequality (\ref{even-1}) is obvious. Otherwise, by
$2^{\frac{\ell-2n+1}{2}}\leq \frac{16}{9}2^{\ell-4}$,
the desired inequality (\ref{even-1}) follows.
If $q\geq3$, then $\sqrt{\Theta}=\sqrt{(2q^{\ell-1}-3q+10)^2+36q-108}
\geq2q^{\ell-1}-3q+10$, which implies that $\frac{1}{9}(2q^{\ell-1}-3q+19)\leq \alpha$,
and therefore the following inequality is sufficient:
\begin{equation}\label{shi-m}
(q-1)^{\frac{n}{2}}q^{\frac{\ell-2n+1}{2}}<\frac{1}{9}(2q^{\ell-1}-3q+19).
\end{equation}
 If $\ell=2n=4$, we can easily prove the inequality (\ref{shi-m}). Otherwise, from $$(q-1)^{\frac{n}{2}}q^{\frac{\ell-2n+1}{2}}<q^{\frac{\ell-n+1}{2}}<2q^{\ell-3}-\frac{1}{3}q\leq\frac{1}{9}(2q^{\ell-1}-3q),$$
we deduce that the desired inequality (\ref{shi-m}) follows.

If $x$ is odd, then the inequality (\ref{shi-g}) is equivalent with
$$h(x,q)=9x^2-(2q^{\ell-1}-3q+34)x+2q^{\ell-1}-5q+33<0.$$
Let
$$\beta=\frac{2q^{\ell-1}-3q+34+\sqrt{\Lambda}}{18}\quad\hbox{with}\quad\Lambda=(2q^{\ell-1}-3q+34)^2-36(2q^{\ell-1}-5q+33).$$
Then $h(\beta,q)=0$. Since $h(2,q)=-2q^{\ell-1}+q+1<0$, the inequality
$(q-1)^{\frac{n}{2}}q^{\frac{\ell-2n+1}{2}}<\beta$
is sufficient.
If $q=2$, then $\sqrt{\Lambda}=\sqrt{2^{2\ell}+20\cdot2^{\ell}-44}> 2^{\ell}-2$, which implies that
$\frac{1}{9}(2^{\ell}+13)<\beta$, and therefore the inequality $2^{\frac{\ell-2n+1}{2}}\leq \frac{1}{9}(2^{\ell}+13)$ is sufficient,
which is a direct corollary of (\ref{even-1}).
If $q=3$, then $\sqrt{\Lambda}=\sqrt{4\cdot3^{2\ell-2}+28\cdot3^{\ell-1}-23}> 2\cdot3^{\ell-1}-1$,
which implies that
$\frac{1}{9}(2\cdot3^{\ell-1}+12)<\beta$, and therefore
the inequality $2^{\frac{n}{2}}3^{\frac{\ell-2n+1}{2}}\leq \frac{1}{9}(2\cdot3^{\ell-1}+12)$ is sufficient,
which is a direct corollary of  (\ref{shi-m}).
If $q\geq4$, then $\sqrt{\Lambda}=\sqrt{(2q^{\ell-1}-3q+16)^2+72q-288}\geq2q^{\ell-1}-3q+16$,
which implies that $\frac{1}{9}(2q^{\ell-1}-3q+25)\leq\beta$,
and therefore the inequality
$(q-1)^{\frac{n}{2}}q^{\frac{\ell-2n+1}{2}}<\frac{1}{9}(2q^{\ell-1}-3q+25)$
is sufficient, which is a direct corollary of  (\ref{shi-m}).

(ii). By (\ref{equa4}), we have
$$\frac{x-1}{\lfloor\frac{3x}{2}\rfloor-2}\Delta(q,3,\ell)-\left(\left\lfloor\frac{3x}{2}\right\rfloor-3\right)s'_2
>xC(q,3,\ell)=xq^{3}(q^{\ell-1}+q^{\ell-2}-1).$$
We claim that
$$xq^{3}(q^{\ell-1}+q^{\ell-2}-1)\geq q^\ell(q^2+q+1)-q(q+1).$$
This inequality is equivalent with
$$x\geq\frac{q^{\ell}(q^2+q+1)-q(q+1)}{q^{\ell+2}+q^{\ell+1}-q^3}.$$
Since $\frac{q^{\ell}(q^2+q+1)-q(q+1)}{q^{\ell+2}+q^{\ell+1}-q^3}< \frac{q^{\ell}(q^2+q+1)}{q^{\ell+2}+q^{\ell+1}-q^\ell}=\frac{q^2+q+1}{q^2+q-1}$, which implies that the  inequality $x\geq\frac{q^2+q+1}{q^2+q-1}$ is sufficient. Since $\frac{q^2+q+1}{q^2+q-1}<2$, the desired result follows.
$\qed$

\begin{lemma}\label{lem4.8}
Let $\ell\geq2n\geq4$ and $2\leq x\leq (q-1)^{\frac{n}{2}}q^{\frac{\ell-2n+1}{2}}$. If ${\cal L}$ is a Cameron-Liebler set
in ${\rm Bil}_q(n,\ell)$ with parameter $x$, then ${\cal L}\supseteq{\cal F}=\{\pi\in{\cal M}_n: \tau\subseteq\pi\}$ for some $\tau\in{\cal M}_1$.
\end{lemma}
\proof Let $\pi\in{\cal L}$ and let $c$ be the maximal number of elements of ${\cal L}$ that are pairwise
disjoint. By Theorem~\ref{lem2.9} (iii),
 the number of subspaces in ${\cal L}$ disjoint to $\pi$  is
$(x-1)\Delta(q,n,\ell).$ Within this collection of subspaces, by Lemma~\ref{lem4.5}, there are at most $c-1\leq\lfloor\frac{3x}{2}\rfloor-2$ subspaces $\pi_1,\ldots,\pi_{c-1}$ that are mutually disjoint. By the pigeonhole principle, we find an $i$ such that $\pi_i$ meets at least $$\frac{x-1}{c-1}\Delta(q,n,\ell)\geq\frac{x-1}{\lfloor\frac{3x}{2}\rfloor-2}\Delta(q,n,\ell)$$ subspaces in ${\cal L}$ that are disjoint to $\pi$. We denote this collection of subspaces disjoint to $\pi$  and meeting $\pi_i$ by ${\cal F}_i$.

Now we claim that ${\cal F}_i$ contains an intersecting family. For every $\pi_j\not=\pi_i$, there are at most $s'_2$ subspaces that meet $\pi_i$ and $\pi_j$. In this way, we find at least
$$\frac{x-1}{\lfloor\frac{3x}{2}\rfloor-2}\Delta(q,n,\ell)-(c-2)s'_2\geq \frac{x-1}{\lfloor\frac{3x}{2}\rfloor-2}\Delta(q,n,\ell)-\left(\left\lfloor\frac{3x}{2}\right\rfloor-3\right)s'_2$$
subspaces in ${\cal L}$ that meet $\pi_i$, are disjoint to $\pi$ and  that are disjoint to $\pi_j$
for all $j\not=i$. We denote this subset of ${\cal F}_i\subseteq{\cal L}$ by ${\cal F}'_i$.
This collection ${\cal F}'_i$ is an intersecting family. In fact, if two subspaces $\sigma_1,\sigma_2$ in ${\cal F}'_i$ would be disjoint, then $(\{\pi_1,\ldots,\pi_{c-1}\}\setminus\{\pi_i\})\cup\{\sigma_1,\sigma_2,\pi\}$ would be a collection of $c+1$ pairwise disjoint subspaces in ${\cal L}$, a contradiction.

By Lemma~\ref{lem4.7},  we have  $\frac{x-1}{\lfloor\frac{3x}{2}\rfloor-2}\Delta(q,n,\ell)-\left(\left\lfloor\frac{3x}{2}\right\rfloor-3\right)s'_2>C(q,n,\ell)+q^{n-1}(q-1)$ if $n\not=3$, and $\frac{x-1}{\lfloor\frac{3x}{2}\rfloor-2}\Delta(q,3,\ell)-(\left\lfloor\frac{3x}{2}\right\rfloor-2)s_2'
>q^\ell(q^2+q+1)-q(q+1)$ if $n=3$. This implies that
$\bigcap_{\pi'\in{\cal F}'_i}\pi'\not=\{0\}$ by Theorem~\ref{lem4.3}.
Let $\tau\in{\cal M}_1$ satisfies $\tau\subseteq\bigcap_{\pi\in{\cal F}'_i}\pi$.
Then ${\cal F}'_i\subseteq{\cal F}$.

We now show that ${\cal L}$ contains ${\cal F}$. If $\sigma\not\in{\cal L}$ satisfies $\tau\subseteq\sigma$, by (\ref{equa4}), $\sigma$  meets at least $\frac{x-1}{\lfloor\frac{3x}{2}\rfloor-2}\Delta(q,n,\ell)-\left(\left\lfloor\frac{3x}{2}\right\rfloor-3\right)s'_2>xC(q,n,\ell)$ subspaces in ${\cal F}'_i$. This contradicts Theorem~\ref{lem2.9} (iii). $\qed$

\begin{thm}\label{lem4.9}
Let $\ell\geq2n\geq4$ and $2\leq x\leq (q-1)^{\frac{n}{2}}q^{\frac{\ell-2n+1}{2}}$.
Then are no non-trivial Cameron-Liebler sets in ${\rm Bil}_q(n,\ell)$ with parameter $x$.
\end{thm}
\proof Suppose that ${\cal L}$ is a Cameron-Liebler set in ${\rm Bil}_q(n,\ell)$ with parameter $x$. We prove that ${\cal L}$ is trivial by using induction. By Lemma~\ref{lem4.8},  ${\cal L}\supseteq{\cal F}=\{\pi\in{\cal M}_n: \tau\subseteq\pi\}$ for some $\tau\in{\cal M}_1$. By Lemma~\ref{lem3.1} (iv), ${\cal L}\setminus{\cal F}$ is a Cameron-Liebler set
in ${\rm Bil}_q(n,\ell)$ with parameter $x-1$. If $x=2$, by Theorem~\ref{lem4.2}, there exists some
$\tau'\in{\cal M}_1$ such that ${\cal L}\setminus{\cal F}={\cal F}'=\{\pi\in{\cal M}_n:\tau'\subseteq\pi\}$,
which implies that ${\cal L}$ is trivial.
If $x\geq 3$,  by the induction hypothesis, ${\cal L}\setminus{\cal F}$ is trivial, which implies that ${\cal L}$ is trivial. $\qed$

\section{Concluding remarks}

In 2010, Wang, Guo and Li \cite{Wang} proved that the set ${\cal M}_m$ has a structure of a symmetric association scheme.
Let
$$
K=\{(i,j): 0\leq i\leq m\wedge(n-m), 0\leq j\leq(m-i)\wedge\ell\},
$$
where $a\wedge b=\min\{a,b\}$.
For $(i,j)\in K$, the relation $R_{(i,j)}$ on
${\cal M}_m$ is defined to be the set of pairs $(\tau,\pi)$ satisfying
$$\dim((\tau+E)/E\cap(\pi+E)/E)=m-i\quad\hbox{and}\quad \dim(\tau\cap\pi)=(m-i)-j.$$
Then $\mathfrak{X}_{m,n,\ell}=({\cal M}_m,\{R_{(i,j)}\}_{(i,j)\in K})$ is a symmetric association scheme and called the {\it association scheme based on the attenuated space} $A_q(n+\ell,E)$. The scheme $\mathfrak{X}_{m,n,\ell}$ is a common generalization of
the Grassmann scheme $\mathfrak{X}_{m,n,0}$ and the bilinear forms scheme $\mathfrak{X}_{n,n,\ell}$. It seems interesting to discuss Cameron-Liebler sets in the scheme $\mathfrak{X}_{m,n,\ell}$.

\section*{Acknowledgments}
The author is indebted to the anonymous reviewers for their detailed reports and constructive suggestions.
The author thanks Professor Ferdinand Ihringer and Professor Alexander L. Gavrilyuk for their various remarks and suggestions
while writing this article. This research is supported by National Natural Science Foundation of China (11971146).

\end{CJK*}

\end{document}